\documentclass{article}

\usepackage[a4paper,top=2.54cm,bottom=2.54cm,left=3.17cm,right=3.17cm,%
            includehead,includefoot]{geometry}

\usepackage{amsmath,amssymb,amsfonts,amsthm}
\usepackage{graphicx}
\usepackage{tabu}
\usepackage{float}
\usepackage{xcolor}
\usepackage[numbers,square,sort&compress]{natbib}
\usepackage{hyperref}
  \hypersetup{colorlinks,citecolor=blue,linkcolor=blue,breaklinks=true}
\usepackage{epstopdf}
\usepackage{yhmath} 

\allowdisplaybreaks

\newtheorem{theorem}{Theorem}[section]

\newtheorem{definition}{Definition}[section]

\newcommand{\bbm}{\begin{bmatrix}}
\newcommand{\ebm}{\end{bmatrix}}

\begin{document}

\title{\uppercase{The K-homology of 2-dimensional Crystallographic Groups}}

\author{
  Hang Wang \thanks{Research Center for Operator Algebras (RCOA) of East China Normal University (ECNU), Shanghai 200241, China. Email: wanghang@math.ecnu.edu.cn}
  \and
  Xiufeng Yao \thanks{Research Center for Operator Algebras (RCOA) of East China Normal University (ECNU), Shanghai 200241, China. Email: yaoze@whut.edu.cn
  }
}

\date{December 7, 2021}

\maketitle

\begin{abstract}
In this paper we compute the topological K-homology of 2-dimensional crystal groups. Our method focuses on the fixed points of group action and simplifies the calculation of the K-homology of universal space. The result also verifies the Baum-Connes Conjecture of 2-dimensional crystal groups.
\end{abstract}

\tableofcontents

\section{Introduction}

A crystal could be thought as a structure consisting of small balls arranged periodically. The balls are constantly arranged in the same way to form the basic pattern (Cell) and spread evenly throughout the space. For example, the crystal structures made of molecules or atoms that are common in physics and chemistry such as diamonds.

The initial motivation for the study of crystal groups came from Hilbert's question 18. Our understanding of the algebraic structure of crystal groups also evolved around Hilbert's 18th question \cite{milnor1976hilbert}. Bieberbach defined the crystal group as a subgroup of the Euclidean group, consisting of the semi-direct product of lattice and point group. He also defined the equivalence of crystal groups by using the conjugate of projective transformations \cite{bieberbach1911bewegungsgruppen}, and proved that there are only finite equivalence classes of crystal groups in n-dimensional space \cite{fricke1912vorlesungen}. This is the first fundamental result of Hilbert's 18th question. And then, Zassenhaus proposed the Jordan-Zassenhaus theorem \cite{zassenhaus1937beweis}, which proves that a geometric equivalence class of crystals contains at most a finite number of arithmetic equivalence classes. This is a fundamental result of integral representation theory and the key to Hilbert's 18th problem.

From the point of view of operator algebra, Zassenhaus's work allows us to take one crystal group as an ordinary finite-dimensional amenable group, that is, a finite locally compact group with Haar measure. B.D. Evans showed that the C*-algebra of a crystal group is a section on vector bundle \cite{evans1982c}. The base space of this vector bundle is the equivalence class of the Pontryagin duality of the lattice about the group actions of the point group , and the fibres are matrices. K. F. Taylor found the way to compute the explicit matrix formulation of the group C*-algebra of any locally compact group with finite-index commutative normal subgroups, and then we could use matrix to describe the C*-algebra of a crystal group \cite{taylor1989c}. J. M. G. Fell and R. S. Doran gave the method to compute the irreducible representations of the group C*-algebra for a crystal group \cite{fell1988representations}. It should be noted that, in Evans's description, the base space can be regarded as the fundamental domain of the orbit space. The group action of point group on a boundary point determines whether its isotropy subgroup is trivial. Eric M. Pohorecky used this method to calculate some examples of 2-dimensional crystal groups \cite{pohorecky1990c}.

Due to the development of string theory and solid state physics, the topological K-theory of crystal groups has received much attention \cite{freed2013twisted}. The difficulty here stems from the special semi-direct product structure of crystal groups, and we usually have to use some cohomology tools to discuss the cases on the boundary of the fundamental domain. Mingze Yang studied the topology of the group C* -algebras of two-dimensional crystal groups and calculated their K- theory \cite{yang1997crossed}. He also found it very difficult in higher dimensional cases. 

To solve this problem, we associate to each crystal group $G$ with a universal space $\underline{E}G$ and the assembly map
\begin{align}\label{bc1}\
  K^{G}_{\bullet}(\underline{E}G) \overset{\mu}{\rightarrow} K_{\bullet}(C_{r}^{*}(G)) .
\end{align}
The left-hand-side is the equivariant K-homology of $\underline{E}G$, and the right-hand-side is the K-theory of the reduced group C*-algebra of $G$. This is what we are familiar with in the formulation of the Baum-Connes Conjecture. 

Using E-theory and KK-theory, Higson and Kasparov proved the Baum-Connes Conjecture for countable discrete groups, which act isometrically and metrically on an infinite-dimensional Euclidean space and include crystal groups \cite{higson1997operator}. Thus we could identify the difficult topolocical K-theory of the group C*-algebra to the better accessible equavariant homology theory given by equivatiant topological K-theory. Davis and Lück calculated the results of both sides of Eq.\ref{bc1} in the case of $Z^{n} \rtimes_{\rho} Z_{p}$( p is a prime). Some physicists have applied similar methods to calculations in crystalline materials \cite{Huang2107}. Recently, Gomi and Thiang introduced T-duality to the study of the Baum-Connes conjecture of crystal groups by using twisted K-theory \cite{gomi2019crystallographic}.

The main structure and motivation of this paper are as follows. Section 2 contains necessary preliminaries about crystal groups and equivariant K-homology. In section 2.1 we review the algebraic structure of crystal groups and relevant cohomology. Section 2.1 includes some theorems about the topological K-homology of finite groups, especially the equations what we need in the computation. Section 3 is the core of this paper. We propose a method that greatly simplifies the calculation of the conjugacy classes of elements of finite order. We propose a method that greatly simplifies the calculation. Basically, this method just needs to find the center of transformation of the plane patterns. We list the detailed calculation results of K-homology for all 17 2-dimensional crystal groups. These results are consistent with the known result about K-theory\cite{yang1997crossed} and verify the Baum-Connes Conjecture. Our original motivation is from the both side of the assembly map \ref{bc1} having the similar form of homology groups. We hope that the calculation will inspire us to explore an explicit description of the assembly map \ref{bc1}.

\section{Preliminaries}

In this sections we recall some relevant definitions and theorems.

\subsection{Crystal Group and Cohomology}

In the previous, the objects arranged periodically include not only the "balls" but also the basic patterns composed of balls, which are generally called "cells". The cells should have special properties that allow them to be arranged in geometric symmetry and periodicity. As the mathematical abstraction of crystal, a crystal group describes the transformations of cells in n-dimensional Euclidean space $\mathbb{R}^{n}$. Furthermore, all these transformations are rigid motions, i.e. isometries in $\mathbb{R}^{n}$. In Euclidean geometry, we have the following assertion \cite{o1966elementary}:

\begin{theorem}
  Let $O(n)$ denote the orthogonal group in $\mathbb{R}^{n}$, $V=\{t(v)| v\in V\}$ denote the vector space of translations to distinguish it from the $\mathbb{R}^{n}$. Every isometry can be written in a unique way as a composition $t(v) \circ \phi$, where $t(v)\in V$ and $\phi\in O(n)$. In addition, all isometries (or rigid motions) of $\mathbb{R}^{n}$ form a group, which we call the Euclidean group $Isom(\mathbb{R}^n)$ and
  \begin{align*}
    (t, \phi)\cdot(t', \phi')&=(t+\phi(t'), \phi \phi'), \\ Isom(\mathbb{R}^{n})&=V \rtimes O(n). 
  \end{align*}
\end{theorem}

A n-dimensional crystal group is the subgroup of $Isom(\mathbb{R}^n)$, Zassenhaus gave the most common definition of crystal group:
\begin{definition}
  \label{defcry}
  An abstract group G is isomorphic to an n-dimensional crystallographic group if and only if G contains a finite index, normal, free abelian subgroup of rank n, that is also maximal abelian.
\end{definition}

Now let us analyze the algebric structure of crystal groups under the framework of Euclidean groups. According to Def.\ref{defcry}, a crystal group fits into an exact sequence
\begin{align}
  M \stackrel{i}{\rightarrow} G \stackrel{p}{\rightarrow} D,
\end{align}
where $i$ and $p$ are natural map. 

The translation part of crystal group generates the free abelian subgroup $M$. This group is often called the lattice. $D$ denotes finite quotient group $G/M$ which is also called the point group, is also a subgroup of $O(n)$. The point group represents rotation and symmetry and generates cells, while the lattice spreads cells all over the space.

Similarly to Euclidean groups, the point group $D$ acts on lattice $M$ by pulling its elements back to G and conjugating. More precisely, if $d\in D$,$a\in M$, we write the action of $d$ on $a$ as
\begin{align}
  d.a =\gamma(d)i(a)\gamma(d)^{-1}=(\varphi(d),d)(a,1)(\varphi(d),d)^{-1}=(d(a),1),
\end{align}
where $\gamma(d)=(\varphi(d),d)$ is a pullback. However, there is a problem we need to point out here. $\gamma$ is given by restricting a pullback from $O(n)$ to $V \rtimes O(n)$ to subgroup $D$. And every pullback is determined by a map $\varphi: D\rightarrow V$. So it may not be a group homomorphism and is not unique. 

The solution to this problem is composing with the exact sequence
\begin{align}\label{eqq1}
  M\rightarrow V \rightarrow V/M.
\end{align}
For an n-dimensional crystal group $G$, considering lattice $M \cong \mathbb{Z}^{n}$ and V being $\mathbb{R}^{n}$, we could write Eq.\ref{eqq1} as
\begin{align}
  \label{eqq2}
  0\rightarrow M\cong \mathbb{Z}^{n} \rightarrow \mathbb{R}^{n}\rightarrow \mathbb{R}^{n} / \mathbb{Z}^{n}\cong \mathbb{T}^{n}\rightarrow 0,
\end{align}
where $\mathbb{T}^{n}$ denotes n-torus. And then we could have a well-defined map $s: D\rightarrow \mathbb{T}^{n}$ associating with $\varphi$. $s$ satifies the condition
\begin{align}
  \label{1cocycle}
  s(d_{1}d_{2})=s(d_{1})+d_{1} (s(d_{2})).
\end{align}
Thus s is a group 1-cocycle. Meanwhile, let $\gamma(d)\gamma(c)=\gamma(dc)\alpha(d,c)$, we could define a 2-cocycle
\begin{align}
  \alpha(d,c):=\gamma(d)\gamma(c) \gamma(dc)^{-1}&=(\varphi(d),d)(\varphi(c),c)(\varphi(dc),dc)^{-1} \\
  &=(\varphi(d)+d(\varphi(c))-\varphi(dc),1).
\end{align}

It is easy to construct cohomology group $H^{\bullet -1}(D,\mathbb{T}^{n})\cong H^{\bullet}(D,\mathbb{Z}^{n}) $ of $G$ with coefficients in $\mathbb{T}^{n}\cong \mathbb{R}^{n}/\mathbb{Z}^{n} \cong V/M$. For a finite group, we often care about its cohomology and homology which are key to K-theory and many other questions. For a crystal group, its cohomology group distinguishes it from other crystal groups with the same lattice and point group. In other words, there exists one-to-one correspondence between crystal groups and the orbits their normalizer acting on the cohomology group $H^{1}(D,\mathbb{T}^{n})$ \cite{schwarzenberger1980n}. We could classify different crystal groups by computing $H^{1}(D,\mathbb{T}^{n})$. Pólya \cite{polya1924xii} and Niggli \cite{NiggliXIIIDF} showed that there are 17 different 2-dimensional crystal groups.

\subsection{Equivarian K-homology of the Crystallographic Groups}
In this part, we will introduce the main theorem we will be using to compute the K-homology of crystal groups. Most definitions and results can be found in \cite{mislin2003proper} and \cite{higson2004group}.

For a discrete group $G$, the Baum-Connes Conjecture allows us to obtain the K-theory of the reduced group C*-algebra $C^{*}_{r}(G)$ of $G$ by computing the equivariant K-homology of the universal space $\underline{E}G$ for proper $G$-actions. First, the definition of proper action is as follows:

\begin{definition}\cite{baum1994classifying}
The $G$-space $X$ is a topolocical space with a given continuous action of $G$ on $X$
$$ G \times X \rightarrow X.$$
The action of $G$ on $X$ is proper if for every $p\in X$, there exists a tripe $(U,H,\rho)$ such that

(i) $U$ is an open neighborhood of $p$, with $gu\in U$for all $(g,u)\in G\times U$,

(ii)$H$ is a compact subgroup of $G$,

(iii)$\rho:U\rightarrow G/H$is a $G$-map from $U$ to the homogeneous space $G/H$.
\end{definition}

The topolocical spaces with group action also have the notion of properness.

\begin{definition}\cite{higson2004group}
A $G$-space $X$ is proper if for every $x\in X $, there is a  $G$-invariant open subset $U\subseteq X$ containing $x$, a finite subgroup $H$ of $G$, and a $G$-equivariant map from $U$to $G/H$.
\end{definition}

Now we can use proper property to define universal spaces. It is important to notice that the universal spaces associated with a group are not unique, and we have group-quivariantly homotopy equivalence. But now we just introduce the notation $\underline{E}G$. Some papers also use universal example to represent it.

\begin{definition}
  The universal space $\underline{E}G$ is a proper $G$-space for proper actions of $G$, if for any proper $G$-space $X$, there exists a $G$-map $f: X\rightarrow\underline{E}G$, and any two $G$-maps from $X$ to $\underline{E}G$ are $G$-homotopic.
\end{definition}

On one hand, \cite{baum1994classifying} shows that an n-dimensional crystal group acts properly on $\mathbb{R}^{n}$ and $\mathbb{R}^{n}$ is a model of $\underline{E}G$ for it. Considering lattice $M\cong \mathbb{Z}^{n}$ acts freely on $\mathbb{R}^{n}\cong\underline{E}G$, the group action of $G$ on $\mathbb{R}^{n}$ could be regarded as the group action of the point group $D$ on torus $\mathbb{T}^{n}\cong \mathbb{R}^{n}/\mathbb{Z}^{n}$ by Eq.\ref{eqq2}. Then we obtain the first main equation in the computation of $K^{G}_{\bullet}(\underline{E}G)$:
\begin{align}\label{Khomology}
  K^{G}_{\bullet}(\underline{E}G)\cong
  K^{\mathbb{Z}^{n}\rtimes D}_{\bullet}(\mathbb{R}^{n})\cong
  K^{D}_{\bullet}(\mathbb{R}^{n}/\mathbb{Z}^{n})\cong
  K^{D}_{\bullet}(\mathbb{T}^{n}).
\end{align}
We will use one of the last three items in computation depending on different crystal groups in the next section.

On the other hand, as a G-CW-complexes, the skeleton filtration of $\underline{E}G$ could be used to construct Bredon homology $H_{i}^{\mathfrak{F i n}}\left(G ; K_{j}^{G}(\cdot)\right)$ which is associated with Chern character. So we could use Chern character and associated equivariant homology group to reduce the computation of $K^{G}_{\bullet}(\underline{E}G)$. The following theorem is originally from \cite{luck2002chern} and summarized by Valette\cite{mislin2003proper}.
\begin{theorem}\cite{mislin2003proper}
  Let $(X,A)$ be a pair of proper G-CW-complexes. Then there is a decomposition
  \begin{align}\label{1}
    C h_{\bullet}^{G}: \bigoplus_{i \in \mathbb{Z}} H_{\bullet+2 i}^{\mathfrak{F i n}}\left(X, A ; \mathbb{Q} \otimes R_{\mathbb{C}}\right) \stackrel{\cong}{\longrightarrow} K_{\bullet}^{G}(X, A) \otimes \mathbb{Q} .
  \end{align}
  For the groups on the left hand side one has isomorphisms
  \begin{align}\label{2}
    H_{k}^{\mathfrak{F} \operatorname{in}}\left(X, A ; \mathbb{Q} \otimes R_{\mathbb{C}}\right) \cong \bigoplus_{[g] \in \operatorname{FC}(G)} H_{k}\left(\left(X^{g}, A^{g}\right) / C_{G}(g) ; \mathbb{Q}\right), \quad k \in \mathbb{N}.
  \end{align}
  $FC(G)$ is the set of conjugacy classes of elements of finite order in $G$. $(\cdot)^{g}$ denotes stabilizer, the set of fixed points of some element $g$. And $C_{G}(g)$ is centralizer of $g$.
\end{theorem}

  For a 2-dimensional crystal group $G=\mathbb{Z}^{2}\rtimes D$, the group action of $G$ on $\underline{E}G\cong \mathbb{R}^{2}$ could decent to $\mathbb{T}^{2}$ as Eq.\ref{Khomology}. So we can rewrite using Eq.\ref{1} and Eq.\ref{2} as
  \begin{align}
    \label{cp}
    K^{G}_{\bullet}(\underline{E}G)\cong \left\{\begin{matrix}
      K^{\mathbb{Z}^{n}\rtimes D}_{\bullet}(\mathbb{R}^{n})\cong \bigoplus_{i \in \mathbb{Z}} \bigoplus_{[g] \in \operatorname{FC}(G)} H_{i+2k}\left(\left(\mathbb{R}^{2}\right)^{g} / C_{G}(g)\right), \\ 
      K^{D}_{\bullet}(\mathbb{T}^{n})\cong \bigoplus_{i \in \mathbb{Z}} \bigoplus_{[g] \in \operatorname{FC}(D)} H_{i+2k}\left(\left(\mathbb{T}^{2}\right)^{g} / C_{D}(g)\right), 
      \end{matrix}\right. \quad k \in \mathbb{N}.
  \end{align}

\section{The K-homology of 2D Crystal Groups}

From the end of the previous section, we know that to compute the K-homology of a crystal group, we need to compute the related homogy group. This requires first to find the conjugacy classes of all elements of finite order in group $G$, and then find the corresponding fixed-point sets and centralizers. The key is to find all conjugacy classes $FC(G)$ or $FC(D)$ in $\mathbb{R}^{2}$ or $\mathbb{T}^{2}$.

\subsection{The Computation of Conjugacy Classes}

For any crystal group $G=A\rtimes D$, let $g=(a,d), a\in A, d\in D$ be an arbitrary element, which means the group action of $G$ is composed by the group action of lattice $A$ and point group $D$. If $G$ is 2-dimensional, the lattice could be writen as $A=<t,s>\cong \mathbb{Z}^{2}$, and $t$, $s$ are linearly independent vectors in $\mathbb{R}^{2}$. Then we could use $\{t,s\}$ as a group of coordinate basis and write an element $ a=mt+ns\in A, (m,n\in\mathbb{Z})$ as $a=(m,n)$ for convenience. Forthermore, $a=(m,n)=((m,n),id_{D})$ denotes the embedding of $a$ in $G$.

For an element $d$ in the point group $D$ with its inverse $d^{-1}$, $\gamma(d)=(u_{d},M_{d})$ denotes the pullback of $d$ from $D$ to $G$. And $u_{d}=(u_{d}^{1},u_{d}^{2})\in \mathbb{R}^{2}$is the translation part, while $M_{d}$ is the transformation matrix in the orthogonal group $O(2)$. From now on, we can write
\begin{align}
  G=A\rtimes D=\{a\cdot d|a\in A,d\in D\}=\{(a,id_{D})\cdot (u_{d},M_{d})=(a+u_{d},M_{d})|a\in A, d\in D \}.
\end{align}
The group action of D on A is 
\begin{align}
  d_{\cdot}a=d^{-1}ad=\gamma(d^{-1})a\gamma(d)=(aM_{d},id_{D}).
\end{align}

let's think for a moment in terms of affine groups. If $y=(y^{1},y^{2}),x=(x^{1},x^{2}),b=(b^{1},b^{2})\in\mathbb{R}^{2}$ is three translation vectors and $(0,0)$ is the coordinate origin, for any martrix $A\in O(2)$, we can define an affine transformation similar to the previous form
\begin{align}
  y=Ax+b=((b^{1},b^{2}),E_{2\times 2})\cdot((0,0),A),
\end{align}
where $E_{2\times 2}$ is $2\times 2$ unitary matrix and could also denotes $id_{D}$. Picking $(m,n)$ as new coordinate origin, we could obtain the new affine translation under new coordidate
\begin{align}\label{1.1}
  y'&=((m,n),E_{2\times 2})((b^{1},b^{2}),E_{2\times 2})((0,0),A)((-m,-n),E_{2\times 2})x \\
  &=((b^{1},b^{2})+(m,n)+ -A(m,n),A)x.
\end{align}

Just like before, using $d$ denotes an orthogonal transformation in a 2-dimensional crystal group $G$, the pullback of $d$ is $\gamma(d)=(u_{d},M_{d})=((u_{d}^{1},u_{d}^{2}),M_{d})$. The group action of lattice $A$ on $d$ gives conjugacy class in the following form:
\begin{align}
  a_{\cdot}d=a^{-1}da=(u_{d}+(-m,-n)+(m,n)M_{d},M_{d}),
\end{align}
where $a=(m,n)=((m,n),id_{D})\in A$. Thus, in analogy to Eq.\ref{1.1}, the group action of a lattice on a point group can be thought of as a corresponding orthogonal transformation at the origin of a point in the lattice.

There are three classes of finite ordered elements of pullback of $D$ in 2-dimensional crystal group $G$, we write $\sigma$, $\rho$ and $d$ respectively. $\sigma$ means rotation, $\rho$ means mirror reflection and $d$ means symmetry about a point. $t=(m,n)\in\mathbb{Z}^{2}$is an arbitrary point in the lattice $A$. Since the order of point group is finite, we can compute conjugacy classes just by conjugating by $t$.

For rotation $\sigma$, $\gamma(\sigma)=((0,0),M_{\sigma})$, we could obtain the conjugacy classes of $\sigma$ in $G$:
\begin{align}
  \begin{aligned}
    t^{-1} \sigma t &=\left((-m,-n), id_{D}\right)\left((0,0), M_{\sigma}\right)\left((m, n), id_{D}\right) \\
    &=\left((-m,-n), id_{D}\right)\left((m, n) M_{\sigma}, M_{\sigma}\right) \\
    &=\left((-m, n)+(m, n)M_{\sigma}, M_{\sigma}\right).
    \end{aligned}
\end{align}
We can easily find the the element of conjugacy class of rotation $\sigma$ could be obtained by changing the rotation origin to $(-m,-n)$.

For mirror reflection $\rho=\left((0,0), M_{\rho}\right)$,
\begin{align}
  t^{-1} \sigma t=\left((-m,-n)+(m, n) M_{\rho}, M_{\rho}\right)=\left((-m, 0)+(m, 0) M_{\rho}, M_{\rho}\right).
\end{align}
This is equivalent to shifting the reflection axis $l_{\rho}$ by $(-m,-n)$.

Finally, for symmetry about the origin point, let $$\gamma(d)=((u,v),M_{d}),\quad M_{d}=\begin{bmatrix}
  -1&0 \\ 
  0&-1 
 \end{bmatrix}.$$
\begin{align}
  \begin{aligned}
    t^{-1} \rho t &=\left((-m,-n), id_{D}\right)\left((u, v), M_{d}\right)\left((m, n), id_{D}\right) \\
    &=\left((-m,-n), id_{D}\right)\left((u, v)+(-m,-n), M_{d}\right) \\
    &=\left(\left(\frac{1}{2} u-m, \frac{1}{2} v-n\right)+\left(m-\frac{1}{2} u, n-\frac{-1}{2} v\right) M_{d}, M_{d}\right).
    \end{aligned}
\end{align}
This is the same as symmetry about $\frac{1}{2}\gamma(d)+(-m,-n)$.

The above calculation about conjugacy classes shows that we only need to find all these three classes of group actions of finite order in one cell. More precisely, if $g=(v,\phi)$ is an element of finite order in $G$( i.e. $(v,\phi)^{k}=(0,id_{D})),k\in\mathbb{N}$), it must be a pullback from $D$ because the orders of translation vectors in lattice are infinite. Since $\phi^{k}=id_{D}$, it satisfies equation
\begin{align}\label{eq}
  v+\phi(v)+\cdots+\phi^{k-1}(v)=0.
\end{align}
So finding all elements of finite order in $G$ just amounts to solve the equation Eq.\ref{eq} for every $\phi \in D$. Considering $\phi(v+\phi(v)+\cdots+\phi^{k-1}(v))=0$, we only need to find all the axes and centers of reflections and centers of rotations in one cell. 

\subsection{K-homology of 2-dimensional Crystal Groups}
In the following, we list what we need to compute the K-homology of 17 2-dimensional crystal groups. We will introduce their point groups and related pullback, and write down the conjugacy classes of elements of finite order. We use the same notation as before for convenience.

According to Eq.\ref{Khomology}, we choose $\mathbb{R}^{2}$ or $\mathbb{T}^{2}$ for convenience. Next, We design a table to list all representative of conjugacy classes, the related stabilizer $X^{g}$, the centralizer $C_{\bullet}(g)$ and their quotient $X^{g}/C_{\bullet}(g)$, and the homology groups of the quotients:$$\bigoplus_{k \in even} \bigoplus_{[g] \in FC[\bullet]} H_{k}\left(X^{g}/C_{\bullet}(g)\right)\quad and\quad \bigoplus_{k \in odd} \bigoplus_{[g] \in FC[\bullet]} H_{k}\left(X^{g}/C_{\bullet}(g)\right).$$
K-homology group $K_{0}(\underline{E}G)$ and $K_{1}(\underline{E}G)$ are given by these direct sums of homology groups. 

In addition, we share some figures to show the patterns of crystal groups with the translation vectors and axes and centre points of rotations and reflections. Let $I^{2}=[-1,1]\times[-1,1]\cong \mathbb{R}^{2}/\mathbb{Z}^{2}\subseteq \mathbb{R}^{2}$. This article also includes the images of $I^{2}$ containing the fundamental domains to help the readers understanding the calculations of homology groups. Unless otherwise specified, all of the points we've given in tables are in $I^{2}$. $I^{2}$ brings us another benefit that the group actions of $D$ on $\mathbb{T}^{2}$ coulb be showed on $I^{2}$ by defining map
\begin{align}
  \begin{aligned}
    \mathbb{T}^{2} &\rightarrow I^{2} \\
    (z,w)=(e^{i\pi t},e^{i\pi s})&\mapsto (t,s).
  \end{aligned}
\end{align}

\subsection*{p1}
\begin{figure}[H]
  \centering
  \begin{minipage}[t]{0.45\textwidth}   
  \centering
  
  \includegraphics[height=5cm]{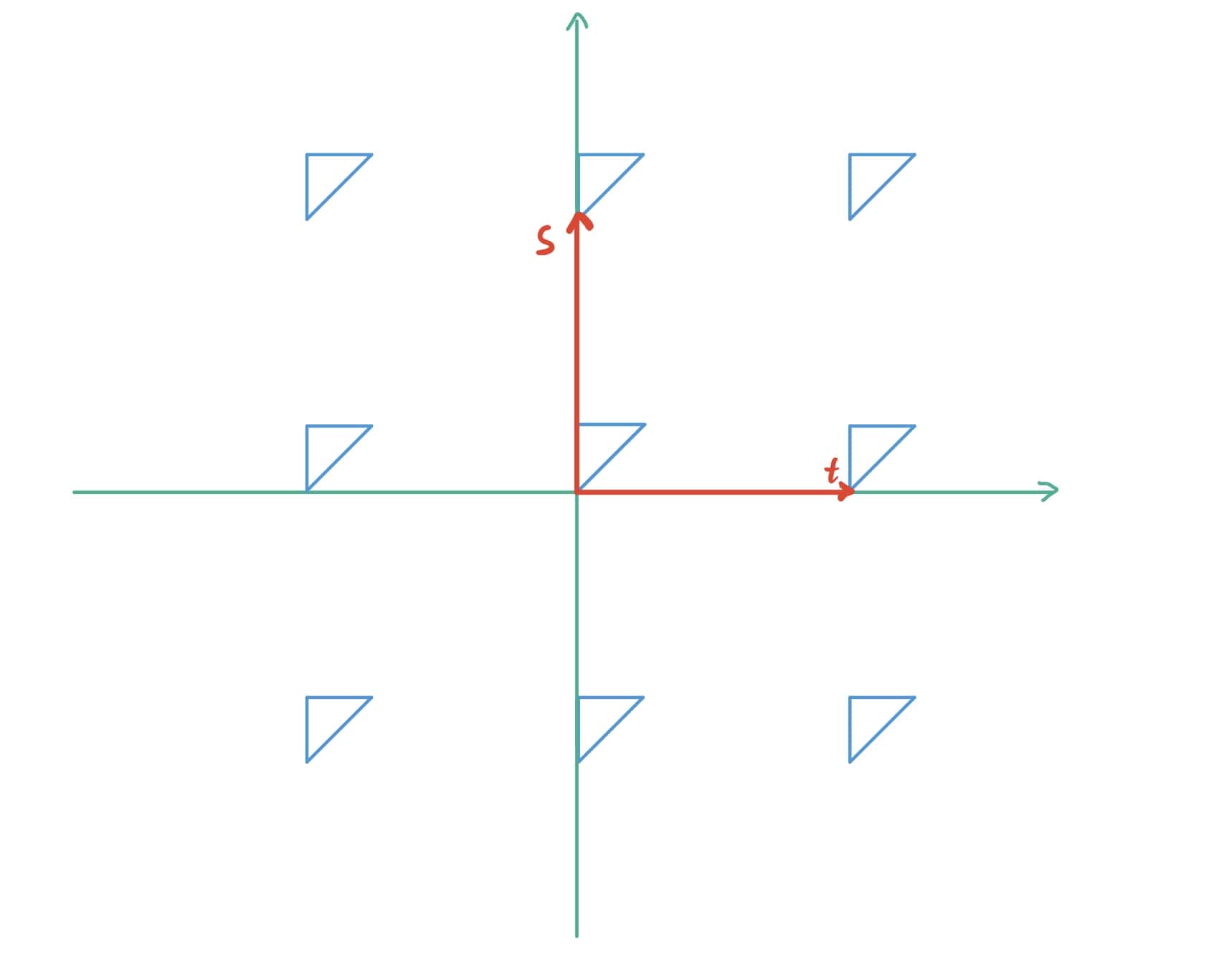}
  
  \caption{p1}
  
  \end{minipage}
  \begin{minipage}[t]{0.45\textwidth}
  
  \centering
  
  \includegraphics[height=5cm]{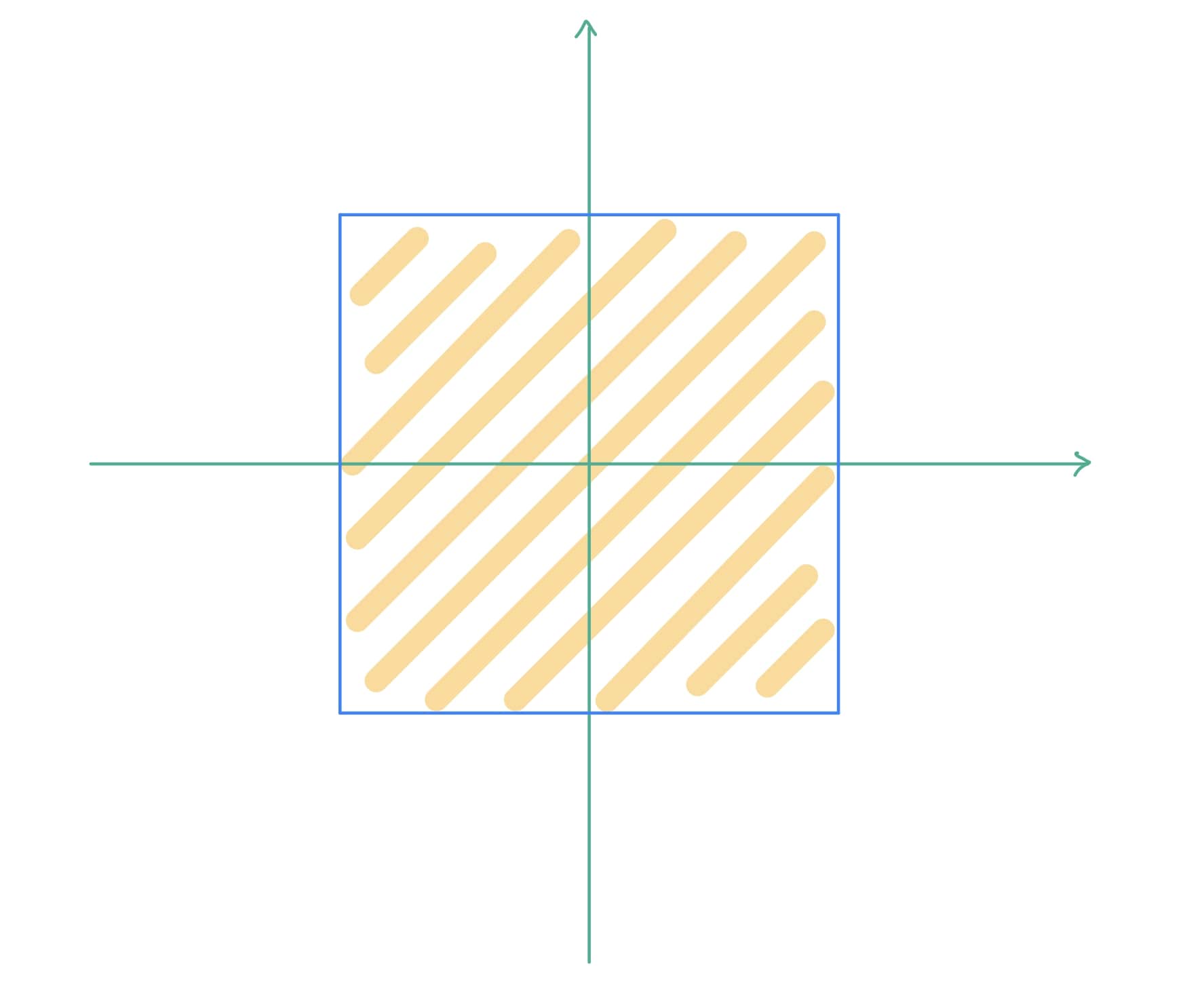}
  
  \caption{Fundamental Domain of p1}
  
  \end{minipage}
  
  \end{figure}

We compute $K^{G}_{\bullet}(\underline{E}G)\cong K^{D}_{\bullet}(\mathbb{T}^{2})$.

1. The group $p1\cong \mathbb{Z}^{2}\cong A=<t,s>$, and its point group $D$ is trivial.
  
2. p1 has no untrivial elements of finite order.

3. See Table \ref{p1}.
\begin{table}[H]
  \footnotesize
  \centering
  \begin{tabu} to 0.95\textwidth{X[c]X[c]}
  \hline
  representative         & identity         \\ \hline
  $X^{g}$                & $\mathbb{T}^{2}$ \\
  $C_{D}(g)$       & $D$          \\
  $X^{g}/C_{D}(g)$ & $\mathbb{T}^{2}$ \\
  \tiny{$\bigoplus_{k \in even} H_{k}\left(X^{g}/C_{D}(g)\right)$} & $\mathbb{Z}\oplus \mathbb{Z}=\mathbb{Z}^{2}$ \\
  \tiny{$\bigoplus_{k \in odd} H_{k}\left(X^{g}/C_{D}(g)\right)$}  & $\mathbb{Z}^{2}$                             \\ \hline
  \end{tabu}%
  \caption{The computation about p1}
  \label{p1}
  \end{table}

  4. $K_{0}(p1)=\mathbb{Z}^{2},K_{1}(p1)=\mathbb{Z}^{2}$.

  \subsection*{p2}
  \begin{figure}[H]
    \centering
    \begin{minipage}[t]{0.45\textwidth}   
    \centering
    
    \includegraphics[height=5cm]{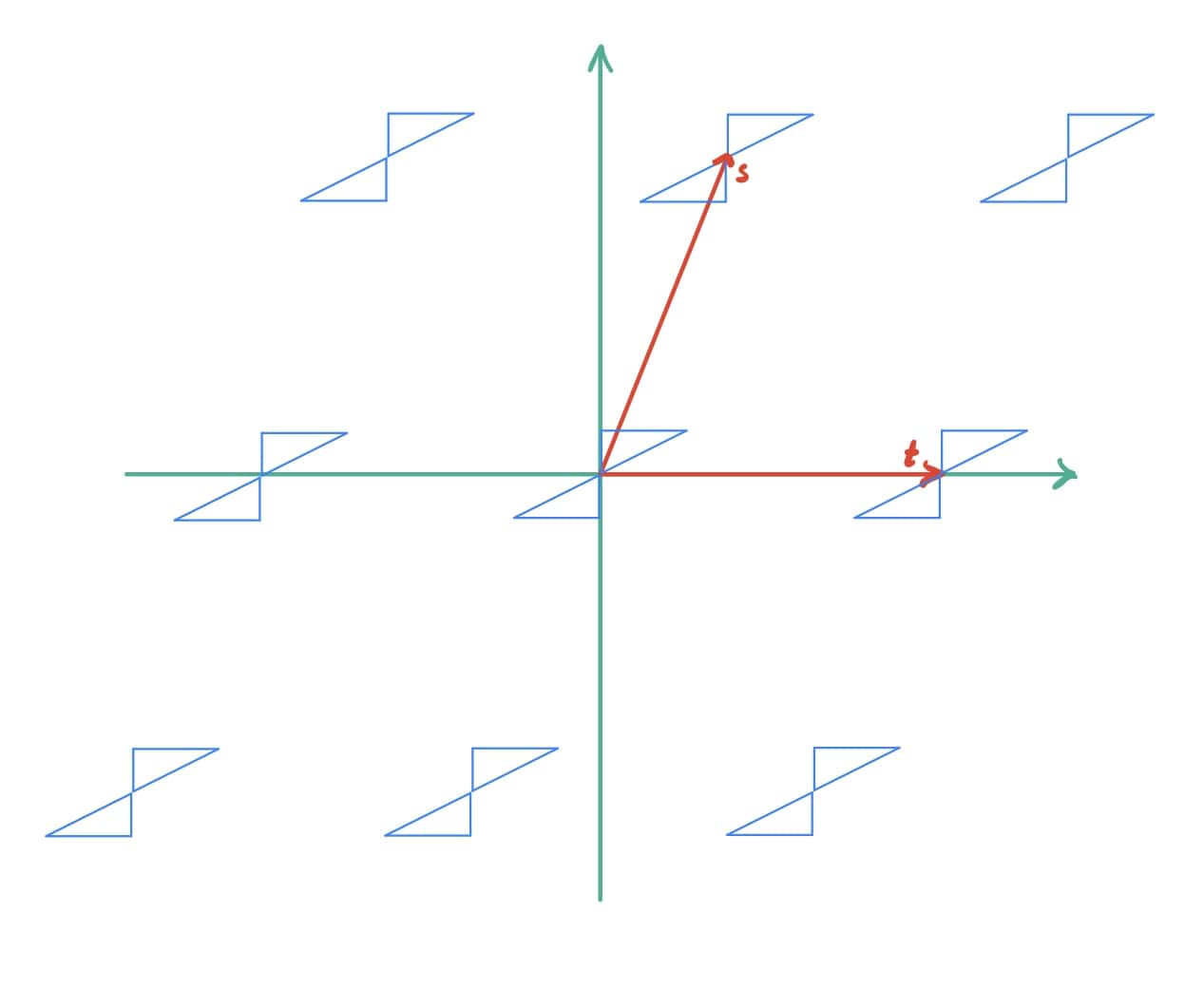}
    
    \caption{p2}
    
    \end{minipage}
    \begin{minipage}[t]{0.45\textwidth}
    
    \centering
    
    \includegraphics[height=5cm]{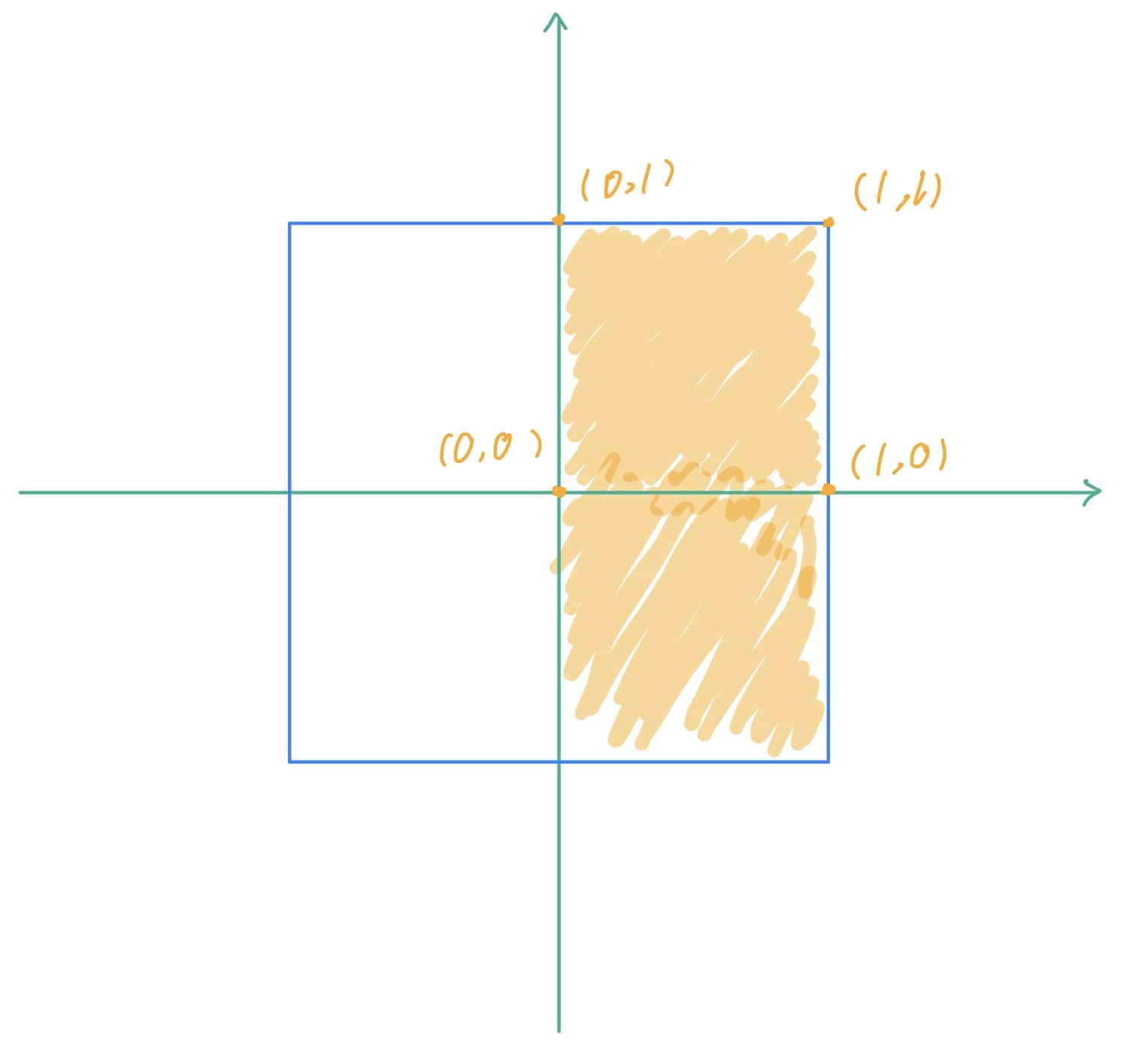}
    
    \caption{Fundamental Domain of p2}
    
    \end{minipage}
    
    \end{figure}
  
  We compute $K^{G}_{\bullet}(\underline{E}G)\cong K^{D}_{\bullet}(\mathbb{T}^{2})$.
  
  1. Group $p2\cong \mathbb{Z}^{2}\rtimes D$. Let $t=(t_{1},0)$ and $s=(s_{1},s_{2})$ be a basis for $A$ ($s_{1},t_{1},t_{2}\neq0$ ). 

  Point group $D=\{id_{D},\sigma\}$ and the pullback:
$$ \gamma(\sigma)=((0,0),M_{\sigma}),\quad M_{\sigma}=\begin{bmatrix}
 -1 & 0\\ 
 0 & -1
 \end{bmatrix}.$$

2. $FC(G)=FC(D)=\{[id_{D}],[\sigma]\}$.
  
  3. See Table \ref{p2}.
  \begin{table}[H]
    \footnotesize
    \centering%
    \begin{tabu} to 0.95\textwidth{X[c]X[c]X[c]}
    \hline
    representative         & identity      & $\sigma$   \\ \hline
    $X^{g}$                & $\mathbb{T}^{2}$ & $(\pm 1,\pm 1)\in \mathbb{T}^{2} $ \\
    $C_{D}(g)$       & $D$    & $D$      \\
    $X^{g}/C_{D}(g)$ & Fundamental Domain & $(1,1),(0,1)$, $(1,0),(0,0) $ \\
    \tiny{$\bigoplus \limits_{k \in even} H_{k}(X^{g}/C_{D}(g))$} & $\mathbb{Z}^{2}$ & $\mathbb{Z}^{4}$ \\
    \tiny{$\bigoplus \limits_{k \in odd} H_{k}(X^{g}/C_{D}(g))$}  & 0 & 0  \\ \hline
    \end{tabu} %
    \caption{The computation about p2}
    \label{p2}
    \end{table}
  
    4. $K_{0}(p2)=\mathbb{Z}^{6}, K_{1}(p2)=0$.

\subsection*{p3}
\begin{figure}[btbp]
  \centering
  \begin{minipage}[t]{0.45\textwidth}   
  \centering
  \includegraphics[height=5cm]{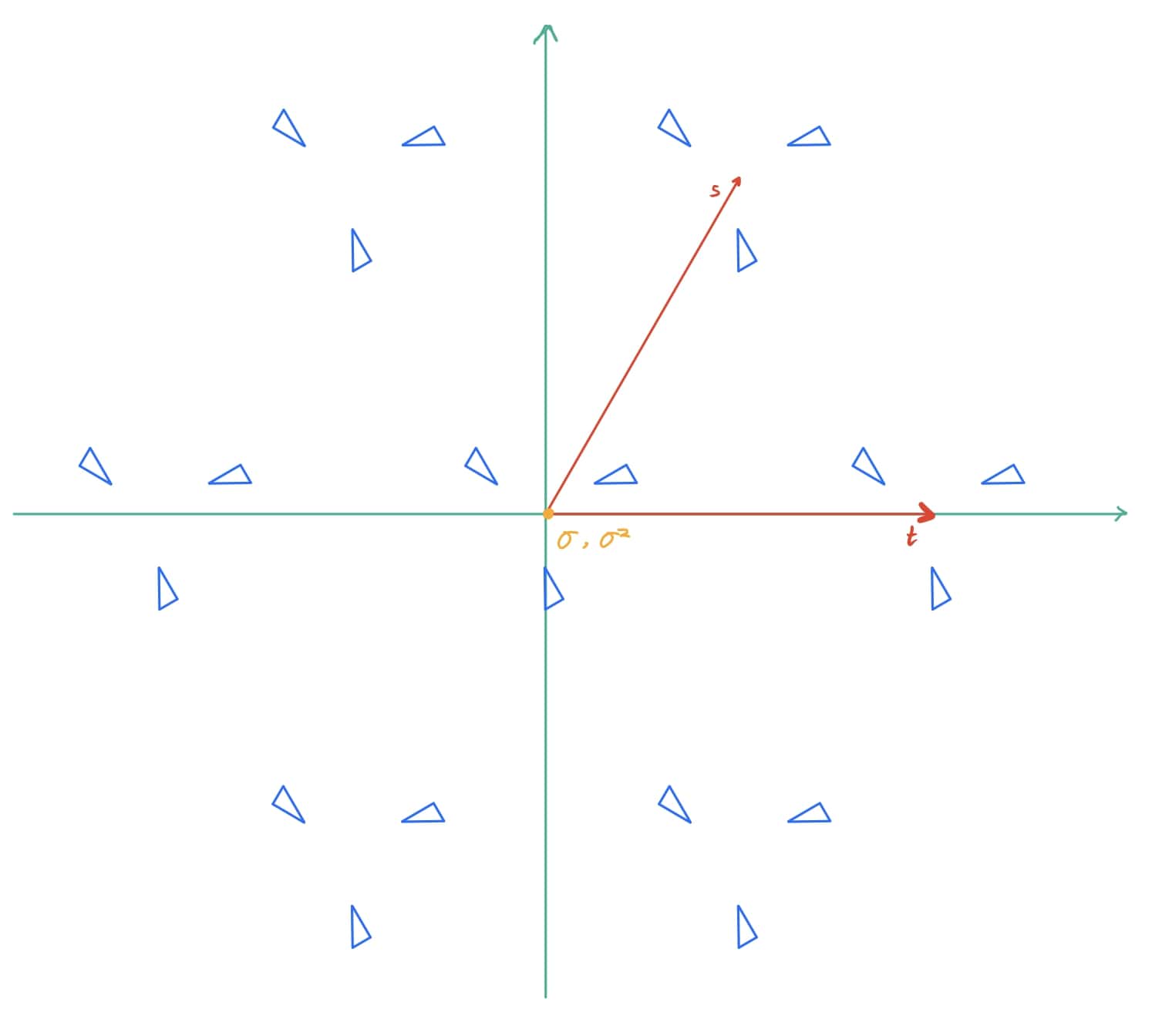}
  \caption{p3}  
  \end{minipage}
  \begin{minipage}[t]{0.45\textwidth}  
  \centering  
  \includegraphics[height=5cm]{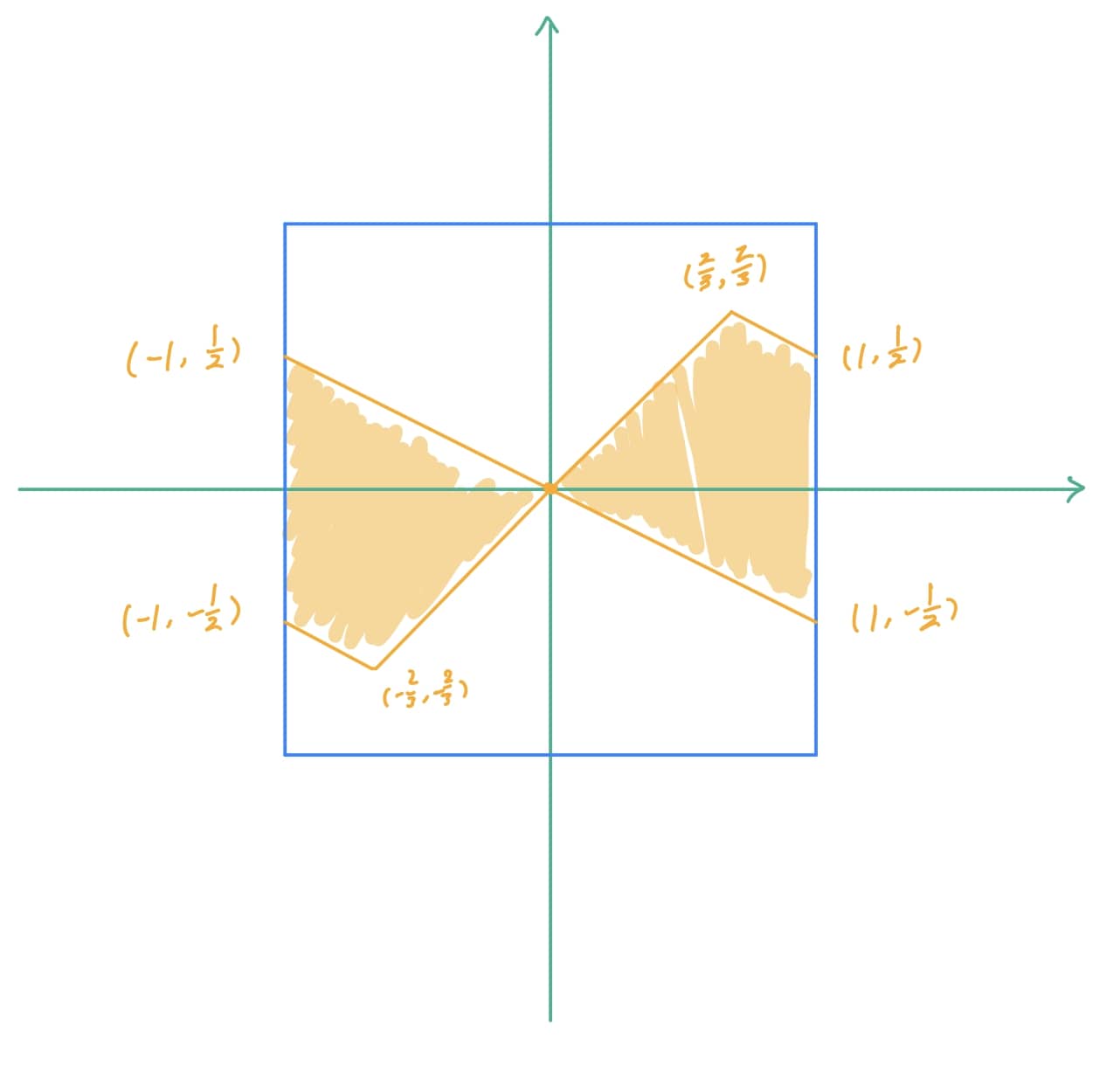}  
  \caption{Fundamental Domain of p3}  
  \end{minipage}  
  \end{figure}

We compute $K^{G}_{\bullet}(\underline{E}G)\cong K^{D}_{\bullet}(\mathbb{T}^{2})$.

1. Let $t=(t_{1},0)$ and $s=(-\frac{1}{2}t_{1},\frac{\sqrt{3}}{2}t_{1})$ form a basis for $A$ ($t_{1}\neq0$ ). 

Point group $D=\{id_{D},\sigma,\sigma^{2}\}$ and the pullback:
\begin{align*}
  \gamma(\sigma)=((0,0),M_{\sigma}),\quad \gamma(\sigma^{2})=((0,0),M_{\sigma^{2}}) \\ M_{\sigma}=\begin{bmatrix}
    0 & 1\\ 
    -1 & -1
    \end{bmatrix}, \quad
    M_{\sigma^{2}}=\begin{bmatrix}
      -1 & -1\\ 
      1 & 0
      \end{bmatrix}.
\end{align*}

2. $FC(G)=FC(D)=\{[id_{D}],[\sigma],[\sigma^{2}]\}$.

3. See Table \ref{p3}.
\begin{table}[H]
  \footnotesize
  \centering%
  \begin{tabu} to 0.95\textwidth{X[c]X[c]X[c]X[c]}
  \hline
  representative         & identity      & $\sigma$ & $\sigma^{2}$  \\ \hline
  $X^{g}$                & $\mathbb{T}^{2}$ & $(\frac{2}{3}, \frac{2}{3})$,$(-\frac{2}{3},-\frac{2}{3})$,$(0,0)$  & same as for $\sigma$\\
  $C_{D}(g)$       & $D$    & $D$      \\
  $X^{g}/C_{D}(g)$ & Fundamental Domain & $(\frac{2}{3}, \frac{2}{3})$,$(-\frac{2}{3},-\frac{2}{3})$,$(0,0)$  & same as for $\sigma$\\
  \tiny{$\bigoplus \limits_{k \in even} H_{k}(X^{g}/C_{D}(g))$} & $\mathbb{Z}^{2}$ & $\mathbb{Z}^{3}$ & $\mathbb{Z}^{3}$ \\
  \tiny{$\bigoplus \limits_{k \in odd} H_{k}(X^{g}/C_{D}(g))$}  & 0 &0 & 0                           \\ \hline
  \end{tabu}%
  \caption{The computation about p3}
  \label{p3}
  \end{table}

  4. $K_{0}(p3)=\mathbb{Z}^{8}, K_{1}(p3)=0$.

  \subsection*{p4}
  \begin{figure}[btbp]
    \centering
    \begin{minipage}[t]{0.45\textwidth}   
    \centering
    \includegraphics[height=5cm]{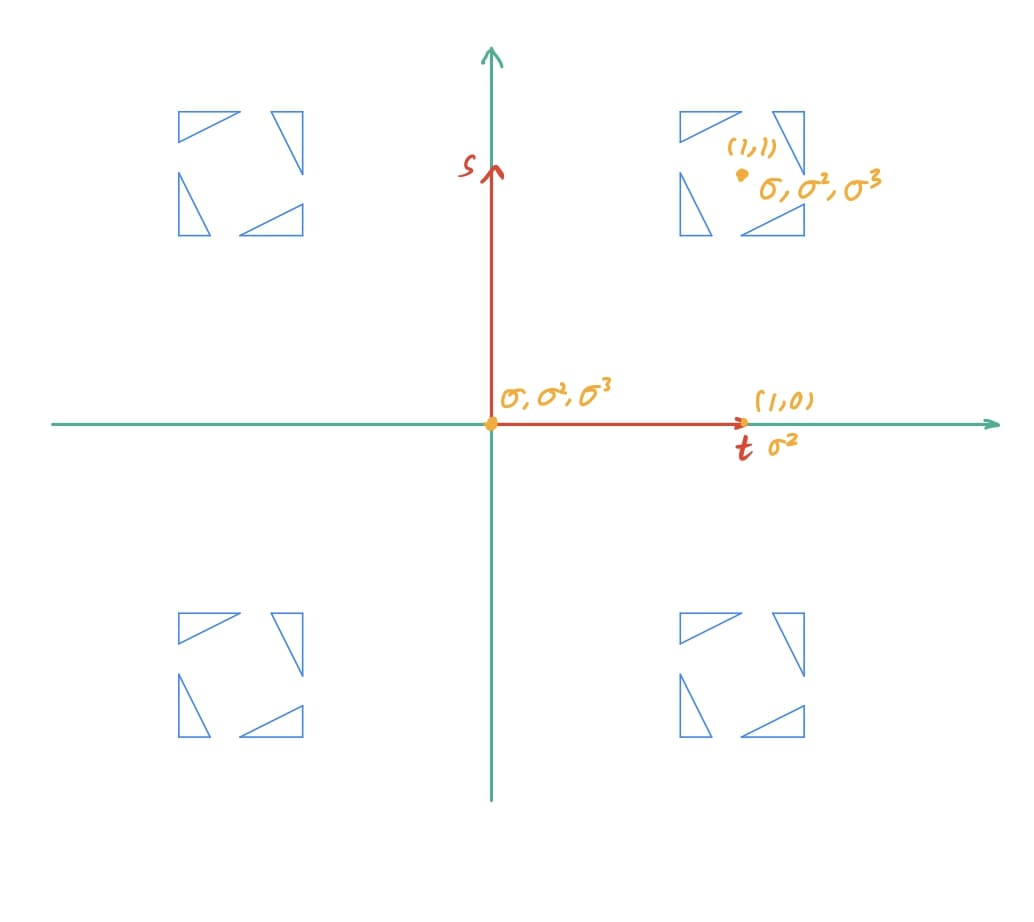}
    \caption{p4}  
    \end{minipage}
    \begin{minipage}[t]{0.45\textwidth}  
    \centering  
    \includegraphics[height=5cm]{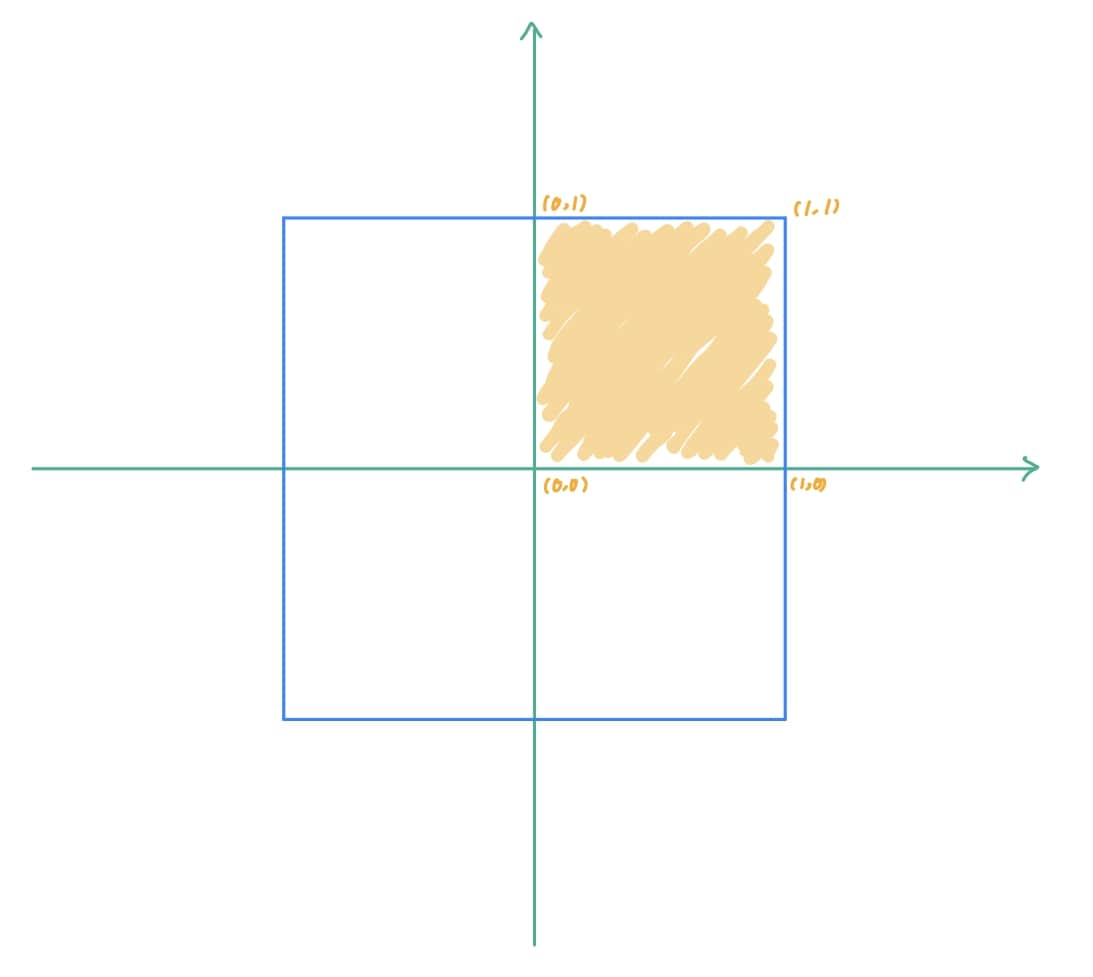}  
    \caption{Fundamental Domain of p4}  
    \end{minipage}  
    \end{figure}

  We compute $K^{G}_{\bullet}(\underline{E}G)\cong K^{D}_{\bullet}(\mathbb{T}^{2})$.

  1. Let $t=(t_{1},0)$ and $s=(0,t_{1})$ form a basis for $A$ ($t_{1}\neq0$ ). 
  
  Point group $D=\{id_{D},\sigma,\sigma^{2},\sigma^{3}\}$ and the pullback:
  \begin{align*}
    \gamma(\sigma)=((0,0),M_{\sigma}),\quad \gamma(\sigma^{2})=((0,0),M_{\sigma^{2}}),\quad \gamma(\sigma^{3})=((0,0),M_{\sigma^{3}}) \\ 
  M_{\sigma}=\begin{bmatrix}
  0 & 1\\ 
  -1 & 0
  \end{bmatrix}, \quad
  M_{\sigma^{2}}=\begin{bmatrix}
    -1 & 0\\ 
    0 & -1
    \end{bmatrix}, \quad
    M_{\sigma^{3}}=\begin{bmatrix}
      0 & -1\\ 
      1 & 0
      \end{bmatrix}.
  \end{align*}
  
  2. $FC(G)=FC(D)=\{[id_{D}],[\sigma],[\sigma^{2}],[\sigma^{3}]\}$.
  
  3. See Table \ref{p4}.
  \begin{table}[H]
    \footnotesize
    \centering%
    \begin{tabu} to 0.95\textwidth{X[c]X[c]X[c]X[c]X[c]}
    \hline
    representative         & identity      & $\sigma$ & $\sigma^{2}$ &$\sigma^{3}$ \\ \hline
    $X^{g}$                & $\mathbb{T}^{2}$ & $(0,0)$, $(1,1)$ & $(1,1),(0,1)$, $(1,0),(0,0)$& $(0,0)$, $(-1,-1)$\\
    $C_{D}(g)$       & $D$    & $D$ &$D$  &$D$   \\
    $X^{g}/C_{D}(g)$ & Fundamental Domain & $(0,0)$, $(1,1)$ & $(1,1),(0,0)$, $(1,0)$& $(0,0)$, $(-1,-1)$\\
    \tiny{$\bigoplus \limits_{k \in even} H_{k}(X^{g}/C_{D}(g))$} & $\mathbb{Z}^{2}$ & $\mathbb{Z}^{2}$ & $\mathbb{Z}^{3}$ &$\mathbb{Z}^{2}$\\
    \tiny{$\bigoplus \limits_{k \in odd} H_{k}(X^{g}/C_{D}(g))$}  & 0 &0 & 0   & 0                        \\ \hline
    \end{tabu}%
    \caption{The computation about p4}
    \label{p4}
    \end{table}
  
    4. $K_{0}(p4)=\mathbb{Z}^{9}, K_{1}(p4)=0$.

    \subsection*{p6}
    \begin{figure}[btbp]
      \centering
      \begin{minipage}[t]{0.45\textwidth}   
      \centering
      \includegraphics[height=5cm]{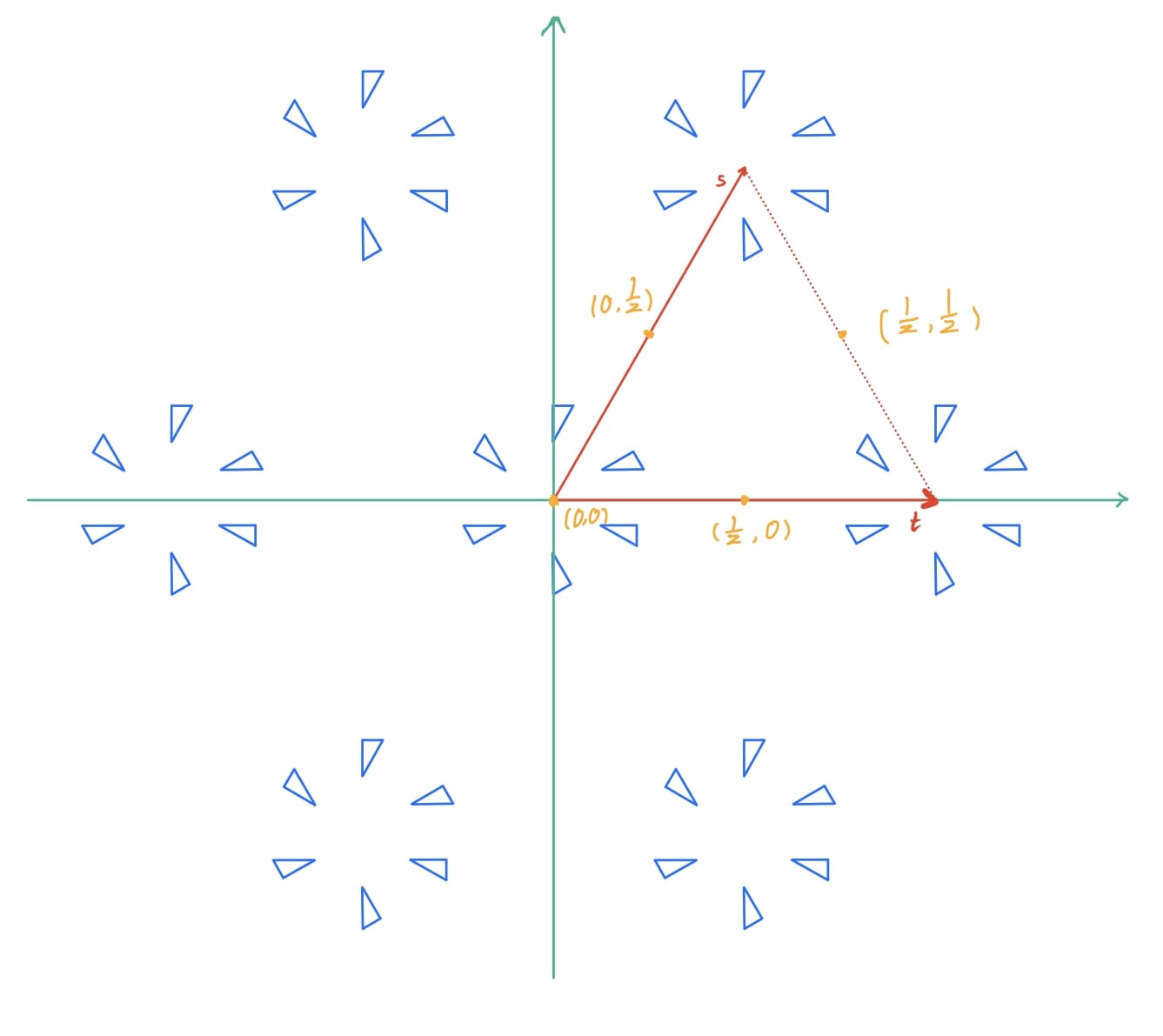}
      \caption{p6}  
      \end{minipage}
      \begin{minipage}[t]{0.45\textwidth}  
      \centering  
      \includegraphics[height=5cm]{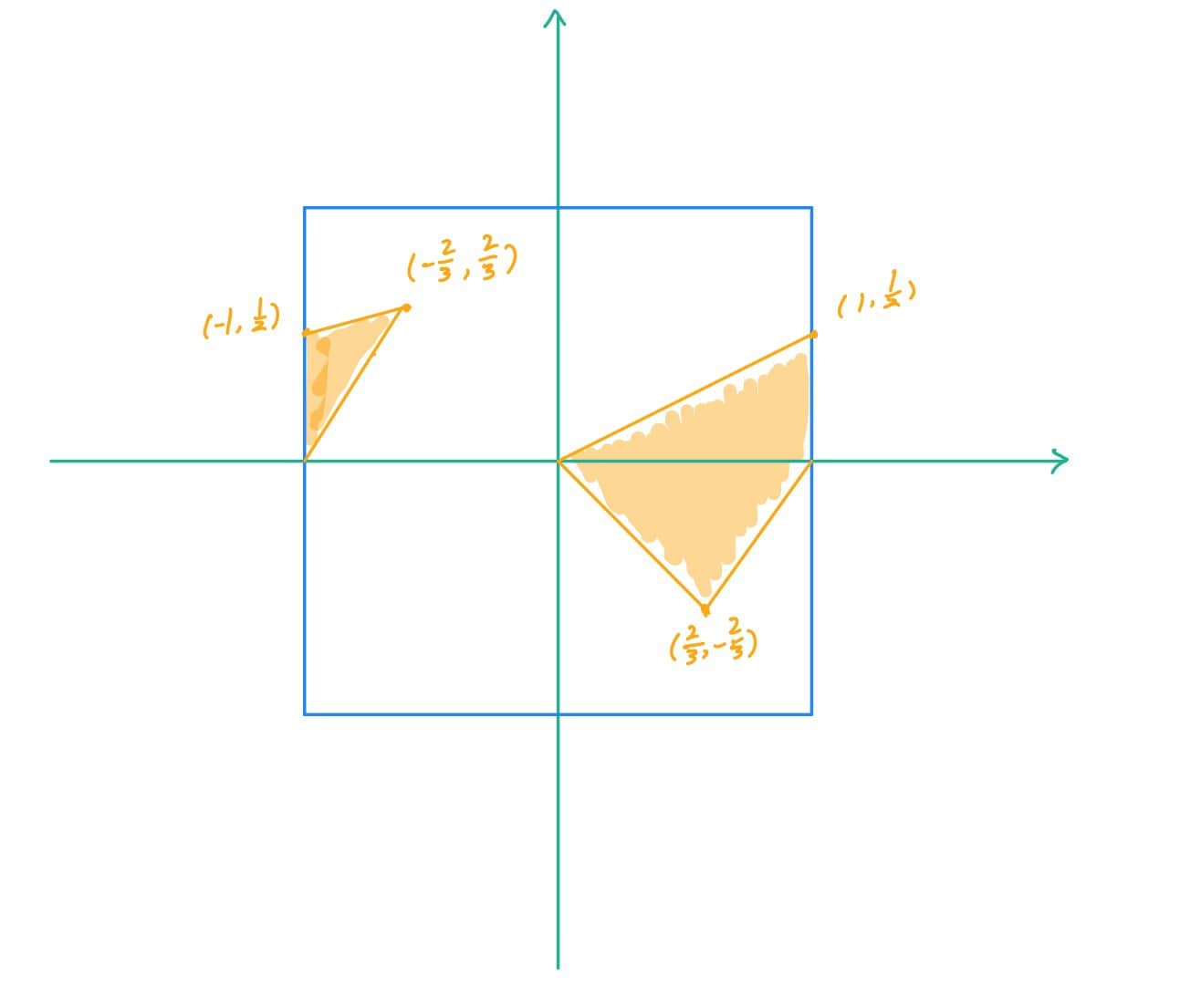}  
      \caption{Fundamental Domain of p6}  
      \end{minipage}  
      \end{figure}
  
    We compute $K^{G}_{\bullet}(\underline{E}G)\cong K^{\mathbb{Z}^{2}\rtimes D}_{\bullet}(\mathbb{R}^{2})$.
  
    1. Let $t=(t_{1},0)$ and  $s=(\frac{1}{2}t_{1},\frac{\sqrt{3}}{2}t_{1})$ form a basis for $A$ ($t_{1}\neq0$ ). 
    
    Point group $D=\{id_{D},\sigma,\sigma^{2},\sigma^{3},\sigma^{4},\sigma^{5}\}$ and the pullback:
    \begin{align*}
      \gamma(\sigma)=((0,0),M_{\sigma}),\quad \gamma(\sigma^{2})=((0,0),M_{\sigma^{2}}),\quad \gamma(\sigma^{3})=((0,0),M_{\sigma^{3}}) \\
      \gamma(\sigma^{4})=((0,0),M_{\sigma^{4}}),\quad \gamma(\sigma^{5})=((0,0),M_{\sigma^{5}})\\
    M_{\sigma}=\begin{bmatrix}
    0 & 1\\ 
    -1 & 1
    \end{bmatrix}, \quad
    M_{\sigma^{2}}=\begin{bmatrix}
      -1 & 1\\ 
      -1 & 0
      \end{bmatrix}, \quad
      M_{\sigma^{3}}=\begin{bmatrix}
        -1 & 0\\ 
        1 & -1
        \end{bmatrix} \\
        M_{\sigma^{4}}=\begin{bmatrix}
          0 & -1\\ 
          1 & -1
          \end{bmatrix}, \quad
          M_{\sigma^{5}}=\begin{bmatrix}
            1 & -1\\ 
            1 & 0
            \end{bmatrix} .
    \end{align*}
    
    2. $FC(G)=\{[id_{D}],[\sigma],[\sigma^{2}],[\sigma^{3}],[\sigma^{4}],[\sigma^{5}],[t\circ \sigma^{3}],[s\circ \sigma^{3}],[s\circ t\circ \sigma^{3}]\}$, where
    $$
\begin{aligned}
t \circ \sigma^{3} &=\left((1,0), M_{\sigma^{3}}\right)=\left(\left(\frac{1}{2}, 0\right)+\left(-\frac{1}{2}, 0\right) M_{\sigma^{3}}, M_{\sigma^{3}}\right) ,\\
s \circ \sigma^{3} &=\left((0,1), M_{\sigma^{3}}\right)=\left(\left(0, \frac{1}{2}\right)+\left(0,-\frac{1}{2}\right) M_{\sigma^{3}}, M_{\sigma^{3}}\right) ,\\
s\circ t\circ \sigma^{3} &=\left((1,1), M_{\sigma^{3}}\right)=\left(\left(\frac{1}{2}, \frac{1}{2}\right)+\left(-\frac{1}{2}, \frac{1}{2}\right) M_{\sigma^{3}}, M_{\sigma^{3}}\right).
\end{aligned}
$$
    
3. See Table \ref{p6}.
 \begin{table}[H]
      \footnotesize
      \centering%
      \begin{tabu} to 0.95\textwidth{X[c]X[c]X[c]X[c]X[c]}
      \hline
      representative         & identity      & $\sigma$ & $\sigma^{2}$ & $\sigma^{3}$  \\ \hline
      $X^{g}$                & $\mathbb{R}^{2}$ & $(0,0)$&$(0,0)$&$(0,0)$\\
      $C_{D}(g)$       & $G$    & $D$ &$D$  &$D$   \\
      $X^{g}/C_{D}(g)$ & Fundamental Domain & $(0,0)$&$(0,0)$&$(0,0)$\\
      \tiny{$\bigoplus \limits_{k \in even} H_{k}(X^{g}/C_{D}(g))$} & $\mathbb{Z}^{2}$ & $\mathbb{Z}$ & $\mathbb{Z}$ &$\mathbb{Z}$ \\
      \tiny{$\bigoplus \limits_{k \in odd} H_{k}(X^{g}/C_{D}(g))$}  & 0 &0 & 0   & 0                        \\ \hline \hline
      $\sigma^{4}$&$\sigma^{5}$&$t\circ \sigma^{3}$&$s\circ \sigma^{3}$&$s\circ t\circ \sigma^{3}$ \\
      $(0,0)$&$(0,0)$&$(\frac{1}{2},0)$&$(0,\frac{1}{2})$&$(\frac{1}{2},\frac{1}{2})$ \\
      $D$&$D$&$id_{D},t\circ \sigma^{3}$&$id_{D},s\circ \sigma^{3}$&$id_{D},s\circ t\circ \sigma^{3}$ \\
      $(0,0)$&$(0,0)$&$(\frac{1}{2},0)$&$(0,\frac{1}{2})$&$(\frac{1}{2},\frac{1}{2})$ \\
      $\mathbb{Z}$ & $\mathbb{Z}$ & $\mathbb{Z}$ &$\mathbb{Z}$&$\mathbb{Z}$ \\
      0&0&0&0&0 \\ \hline
      \end{tabu}%
      \caption{The computation about p6}
      \label{p6}
      \end{table}
    
      4. $K_{0}(p6)=\mathbb{Z}^{10}, K_{1}(p6)=0$.


      \subsection*{pm}
      \begin{figure}[btbp]
        \centering
        \begin{minipage}[t]{0.45\textwidth}   
        \centering
        \includegraphics[height=5cm]{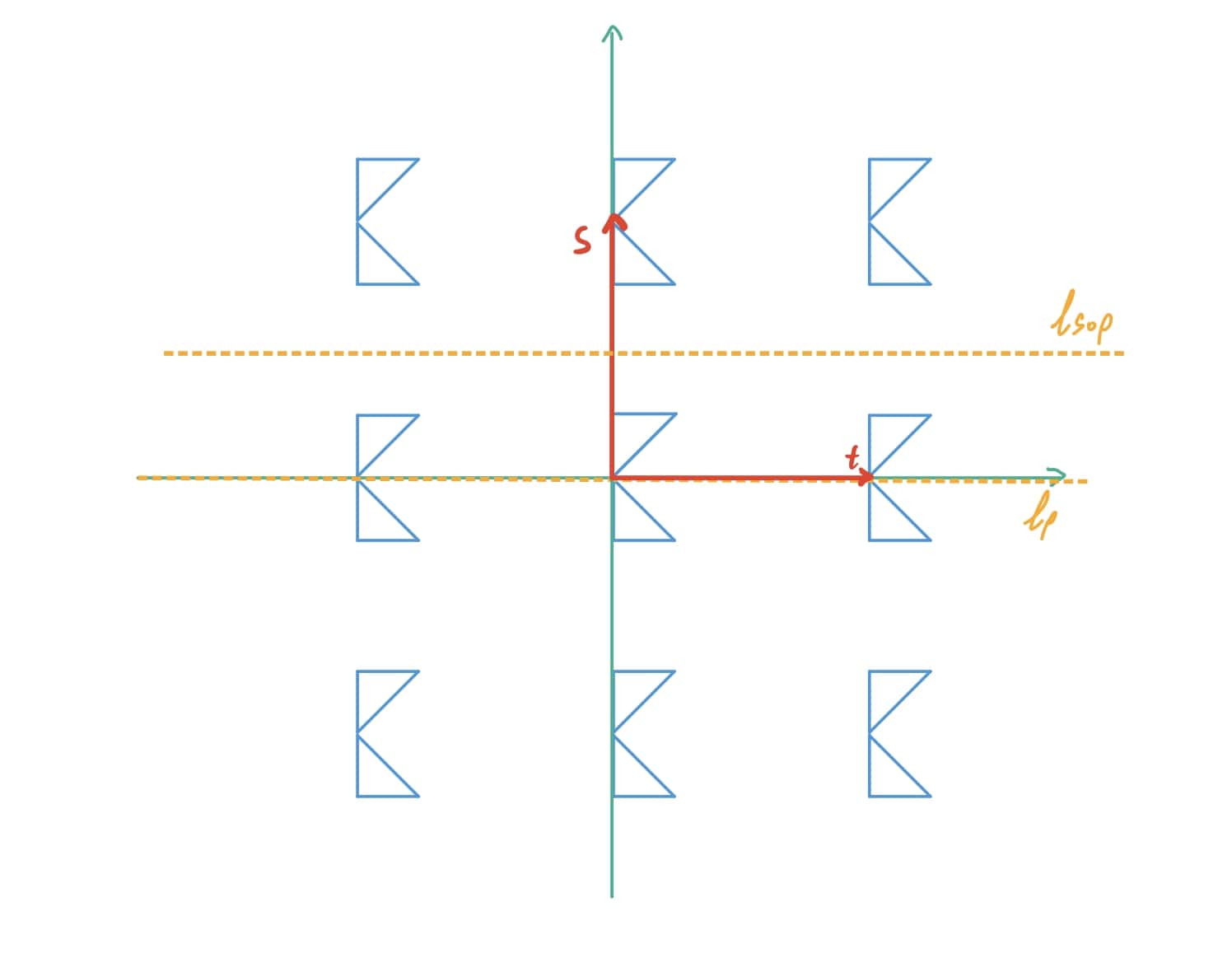}
        \caption{pm}  
        \end{minipage}
        \begin{minipage}[t]{0.45\textwidth}  
        \centering  
        \includegraphics[height=5cm]{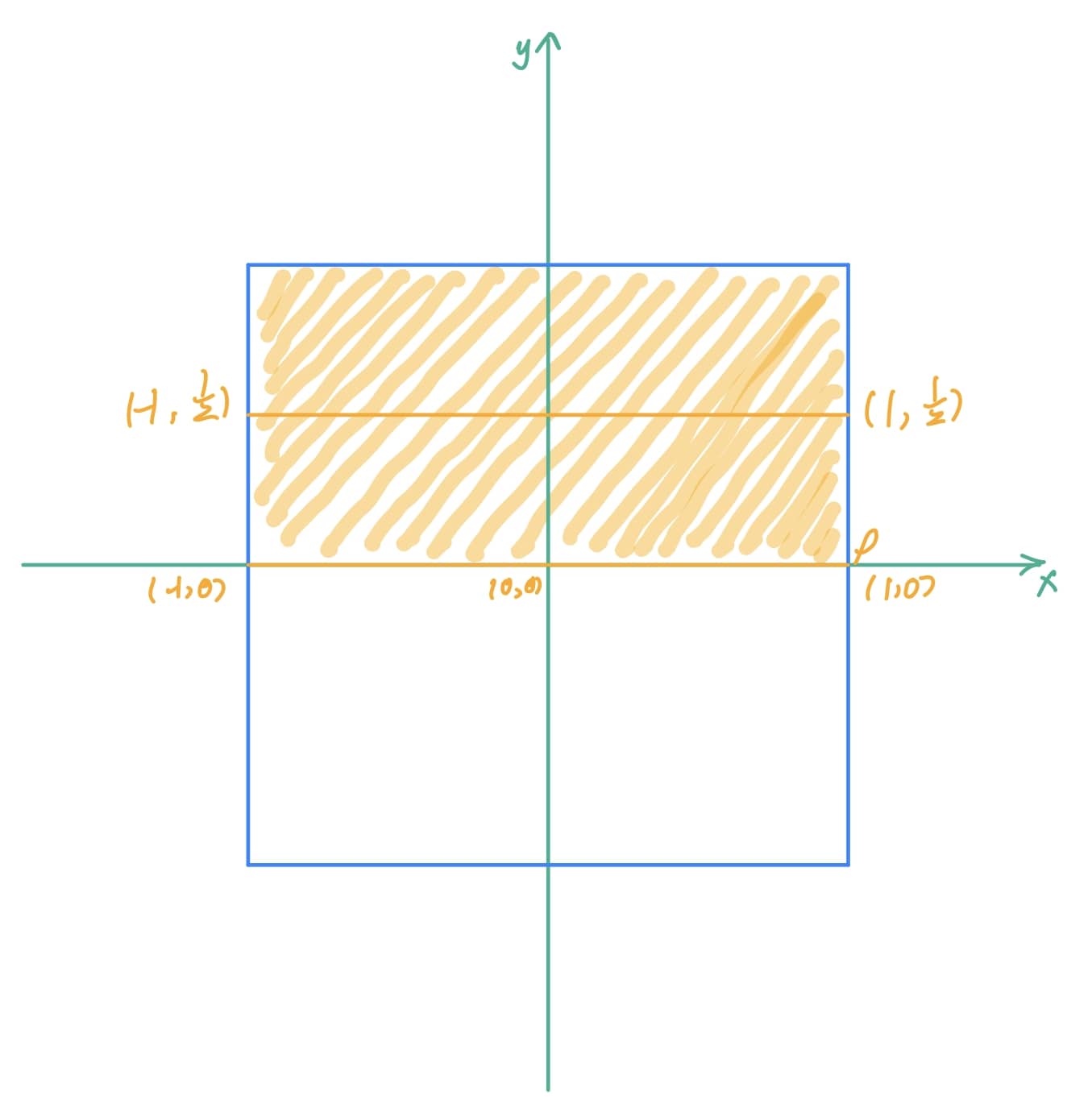}  
        \caption{Fundamental Domain of pm}  
        \end{minipage}  
        \end{figure}
    
        We compute $K^{G}_{\bullet}(\underline{E}G)\cong K^{\mathbb{Z}^{2}\rtimes D}_{\bullet}(\mathbb{R}^{2})$.
    
      1. Let $t=(t_{1},0)$ and $s=(0,s_{2})$ form a basis for $A$ ($t_{1},s_{2}\neq0$ ). 
      
      Point group $D=\{id_{D},\rho\}$ and the pullback:
      \begin{align*}
        \gamma(\rho)=((0,0),M_{\rho}),\quad
      M_{\rho}=\begin{bmatrix}
      1 & 0\\ 
      0 & -1
      \end{bmatrix}.
      \end{align*}
      
      2. $FC(G)=\{[id_{D}],[\rho],[s\circ \rho]\}$.
      
      3. See Table \ref{pm}.
      \begin{table}[H]
        \footnotesize
        \centering%
        \begin{tabu} to 0.95\textwidth{X[c]X[c]X[c]X[c]}
        \hline
        representative         & identity      & $\rho$ & $s\circ \rho$ \\ \hline
        $X^{g}$                & $\mathbb{R}^{2}$ & $(x,0)$& $(x,-\frac{1}{2})$\\
        $C_{D}(g)$       & $G$    & $t,id_{D},\rho$ &$t,\rho$   \\
        $X^{g}/C_{D}(g)$ & Fundamental Domain & $S^{1}$ & $S^{1}$ \\
        \tiny{$\bigoplus \limits_{k \in even} H_{k}(X^{g}/C_{D}(g))$}  & $\mathbb{Z}$ & $\mathbb{Z}$ &$\mathbb{Z}$\\
        \tiny{$\bigoplus \limits_{k \in odd} H_{k}(X^{g}/C_{D}(g))$}  & $\mathbb{Z}$ & $\mathbb{Z}$ &$\mathbb{Z}$                     \\ \hline
        \end{tabu}%
        \caption{The computation about pm}
        \label{pm}         
        \end{table}
      
        4. $K_{0}(pm)=\mathbb{Z}^{3}, K_{1}(pm)=\mathbb{Z}^{3}$.

        \subsection*{cm}
        \begin{figure}[btbp]
          \centering
          \begin{minipage}[t]{0.45\textwidth}   
          \centering
          \includegraphics[height=5cm]{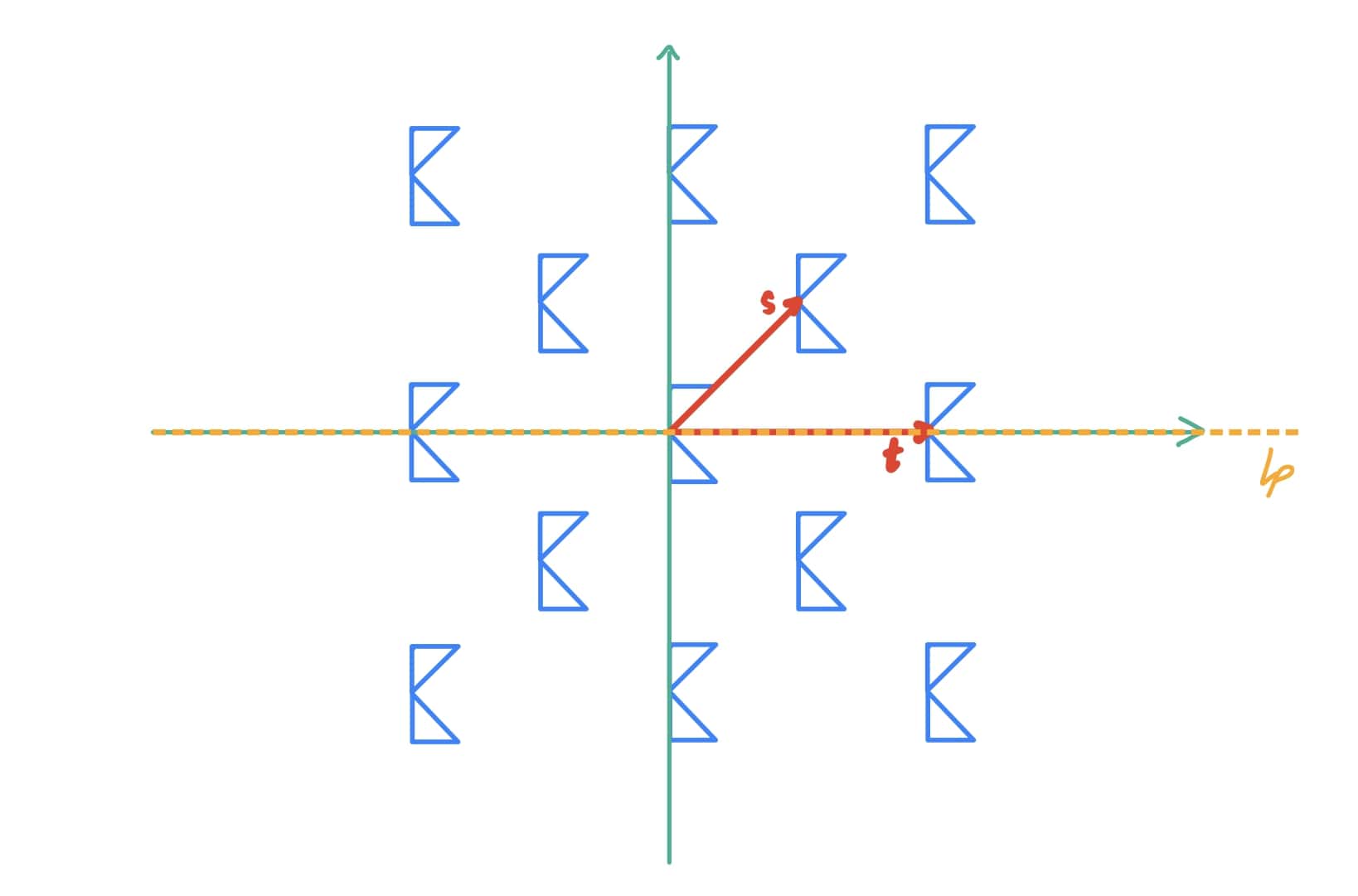}
          \caption{cm}  
          \end{minipage}
          \begin{minipage}[t]{0.45\textwidth}  
          \centering  
          \includegraphics[height=5cm]{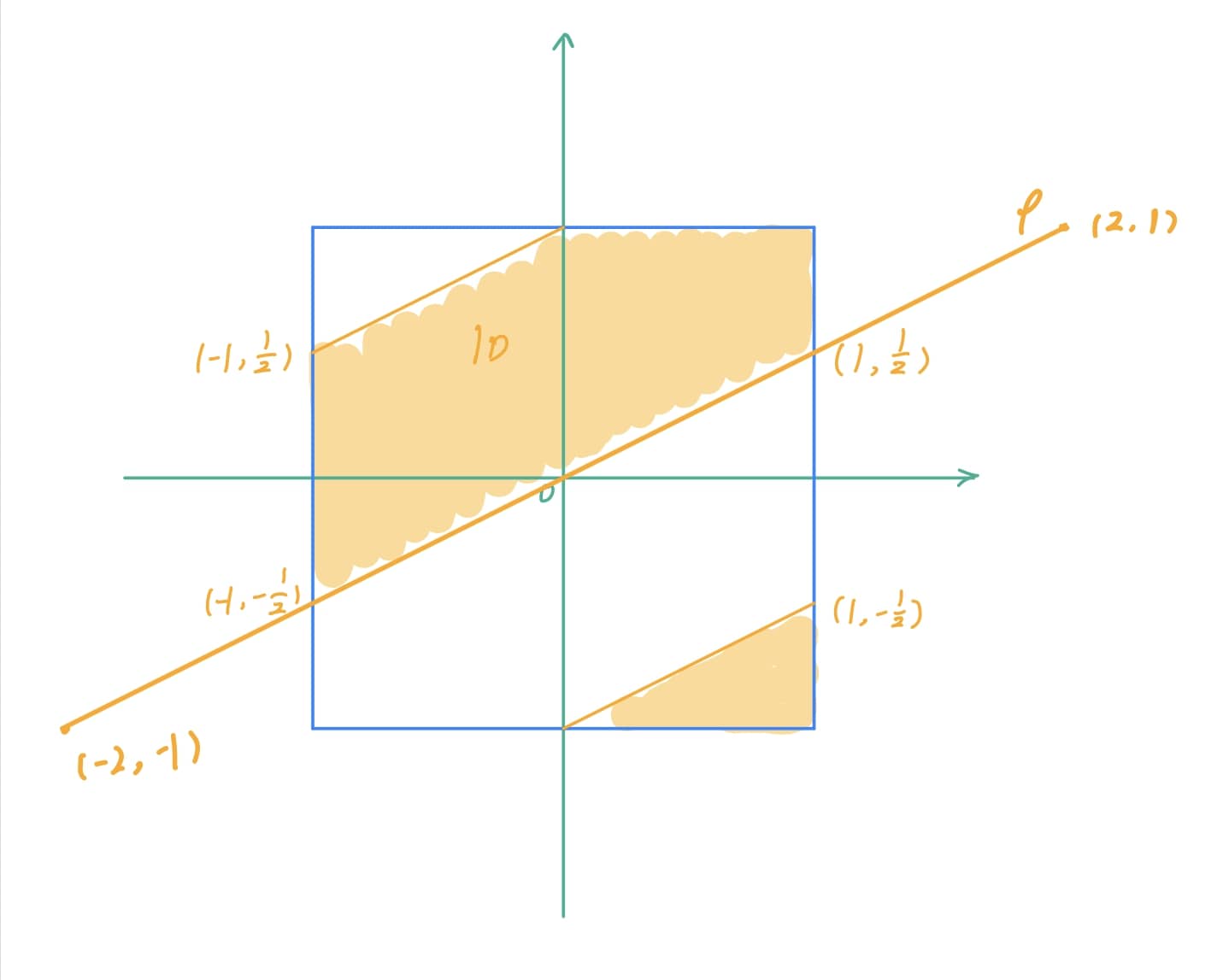}  
          \caption{Fundamental Domain of cm}  
          \end{minipage}  
          \end{figure}
      
          We compute $K^{G}_{\bullet}(\underline{E}G)\cong K^{\mathbb{Z}^{2}\rtimes D}_{\bullet}(\mathbb{R}^{2})$.
      
        1. Let $t=(t_{1},0)$ and $s=(\frac{1}{2}t_{1},\frac{1}{2}s_{2})$ form a basis for $A$ ($t_{1},s_{2}\neq0$ ).
        
        Point group $D=\{id_{D},\rho\}$ and the pullback:
        \begin{align*}
          \gamma(\rho)=((0,0),M_{\rho}),\quad
        M_{\rho}=\begin{bmatrix}
        1 & 0\\ 
        1 & -1
        \end{bmatrix}.
        \end{align*}
        
        2. $FC(G)=\{[id_{D}],[\rho]\}$.
        
        3. See Table \ref{cm}.
        \begin{table}[H]
          \footnotesize
          \centering%
          \begin{tabu} to 0.95\textwidth{X[c]X[c]X[c]}
          \hline
          representative         & identity      & $\rho$  \\ \hline
          $X^{g}$                & $\mathbb{R}^{2}$ & $(x,\frac{x}{2})$\\
          $C_{D}(g)$       & $G$    & $2t+s,id_{D},\rho$    \\
          $X^{g}/C_{D}(g)$ & Fundamental Domain & $S^{1}$  \\
          \tiny{$\bigoplus \limits_{k \in even} H_{k}(X^{g}/C_{D}(g))$}  & $\mathbb{Z}$ & $\mathbb{Z}$ \\
          \tiny{$\bigoplus \limits_{k \in odd} H_{k}(X^{g}/C_{D}(g))$}  & $\mathbb{Z}$ & $\mathbb{Z}$      \\ \hline
          \end{tabu}%
          \caption{The computation about cm}
          \label{cm}         
          \end{table}
        
          4. $K_{0}(cm)=\mathbb{Z}^{2}, K_{1}(cm)=\mathbb{Z}^{3}$.

          \subsection*{pg}
          \begin{figure}[btbp]
            \centering
            \begin{minipage}[t]{0.45\textwidth}   
            \centering
            \includegraphics[height=5cm]{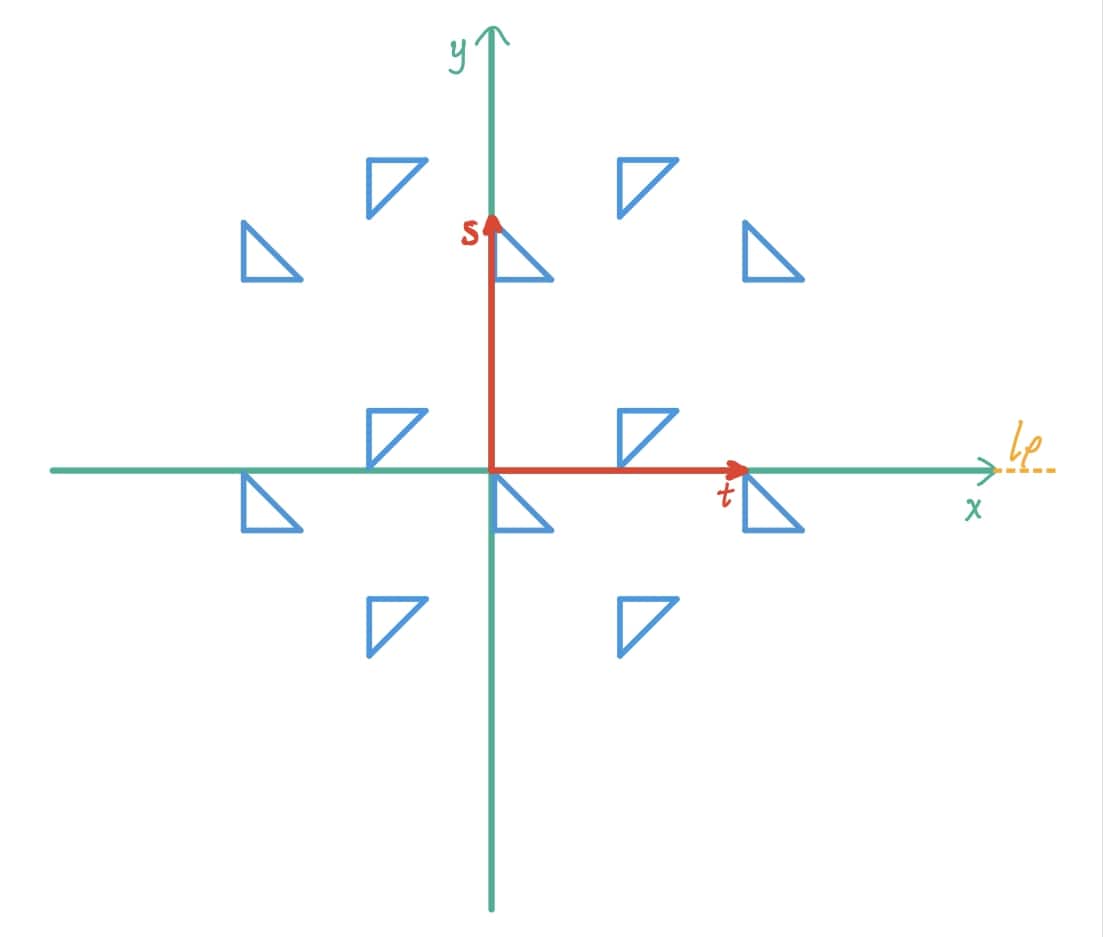}
            \caption{pg}  
            \end{minipage}
            \begin{minipage}[t]{0.45\textwidth}  
            \centering  
            \includegraphics[height=5cm]{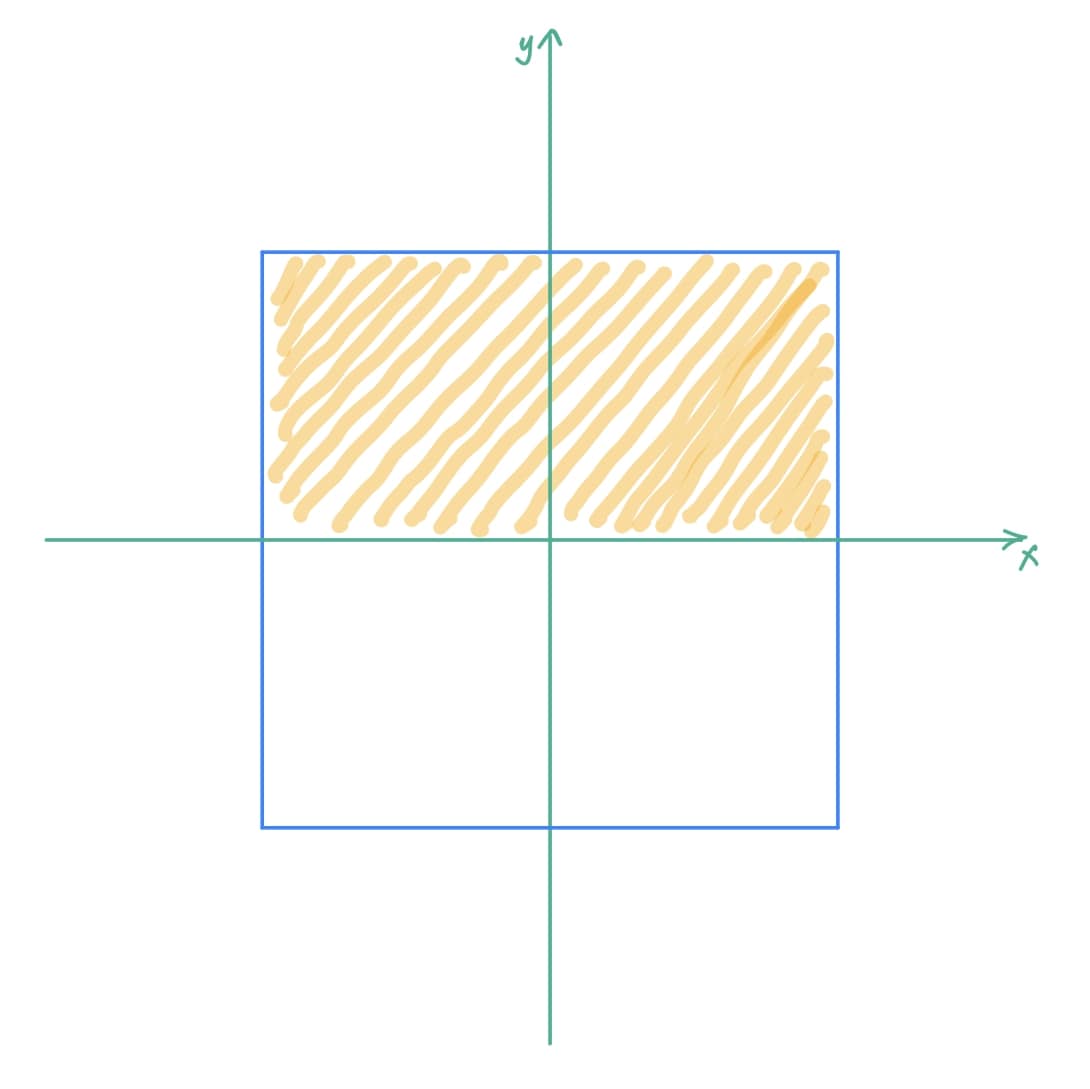}  
            \caption{Fundamental Domain of pg}  
            \end{minipage}  
            \end{figure}
        
            We compute $K^{G}_{\bullet}(\underline{E}G)\cong K^{\mathbb{Z}^{2}\rtimes D}_{\bullet}(\mathbb{R}^{2})$.
        
          1. Let $t=(t_{1},0)$ and $s=(0,s_{2})$ form a basis for $A$ ($t_{1},s_{2}\neq0$ ).
          
          Point group $D=\{id_{D},\rho\}$ and the pullback:
          \begin{align*}
            \gamma(\rho)=((\frac{1}{2},0),M_{\rho}),\quad
          M_{\rho}=\begin{bmatrix}
          1 & 0\\ 
          o & -1
          \end{bmatrix}.
          \end{align*}
          
          2. $FC(G)=\{[id_{D}]\}$.
          
          3. See Table \ref{pg}.
          \begin{table}[H]
            \footnotesize
            \centering%
            \begin{tabu} to 0.95\textwidth{X[c]X[c]}
            \hline
            representative         & identity       \\ \hline
            $X^{g}$                & $\mathbb{R}^{2}$ \\
            $C_{D}(g)$       & $G$      \\
            $X^{g}/C_{D}(g)$ & Cylindar   \\
            \tiny{$\bigoplus \limits_{k \in even} H_{k}(X^{g}/C_{D}(g))$}  & $\mathbb{Z}$ \\
            \tiny{$\bigoplus \limits_{k \in odd} H_{k}(X^{g}/C_{D}(g))$}  & $\mathbb{Z}\oplus \mathbb{Z}_{2} $    \\ \hline
            \end{tabu}%
            \caption{The computation about pg}
            \label{pg}         
            \end{table}
          
            4. $K_{0}(pg)=\mathbb{Z}, K_{pg}(cm)=\mathbb{Z}\oplus \mathbb{Z}_{2}$.

      \subsection*{cmm2}
    \begin{figure}[btbp]
      \centering
      \begin{minipage}[t]{0.45\textwidth}   
      \centering
      \includegraphics[height=5cm]{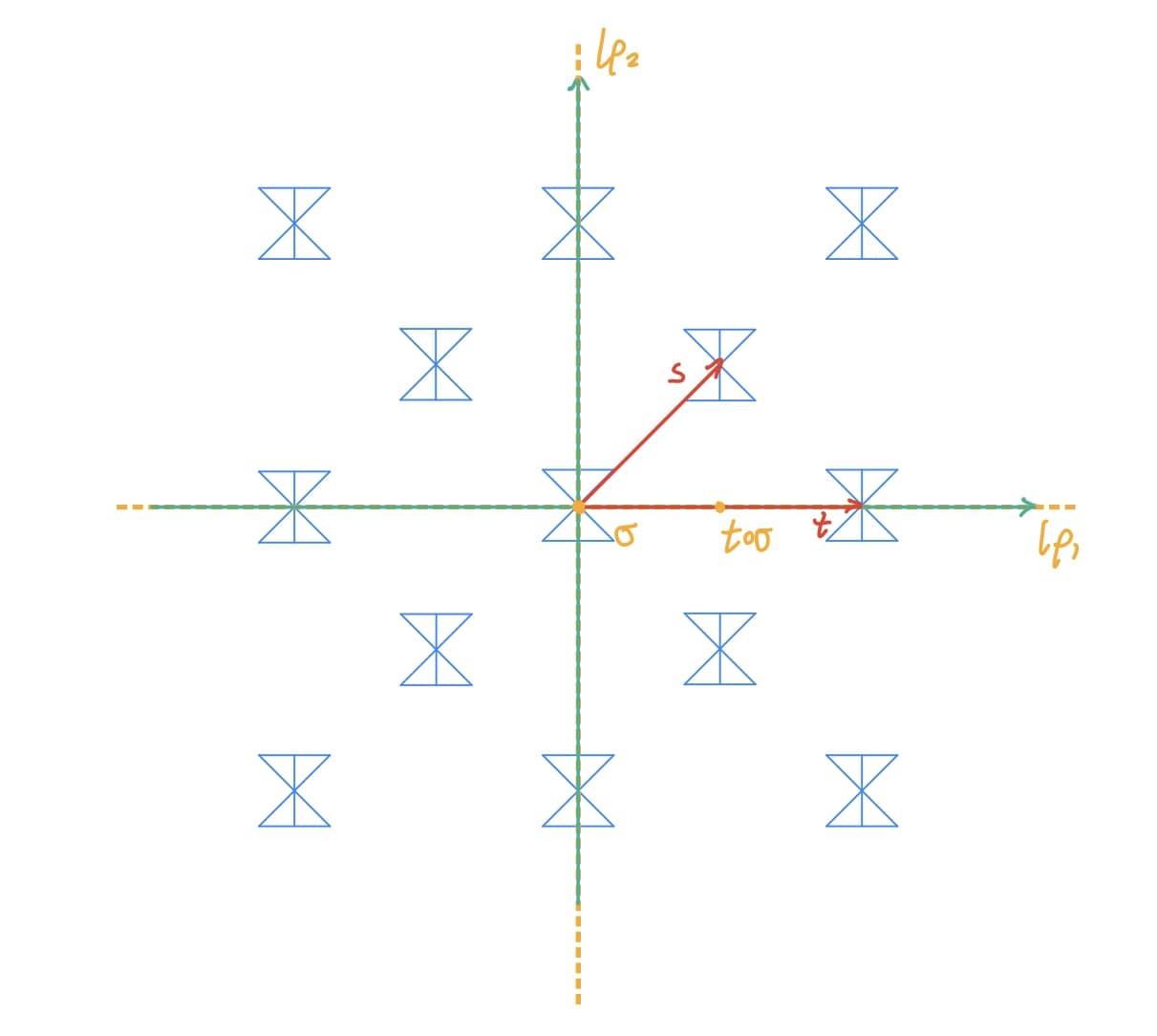}
      \caption{cmm2}  
      \end{minipage}
      \begin{minipage}[t]{0.45\textwidth}  
      \centering  
      \includegraphics[height=5cm]{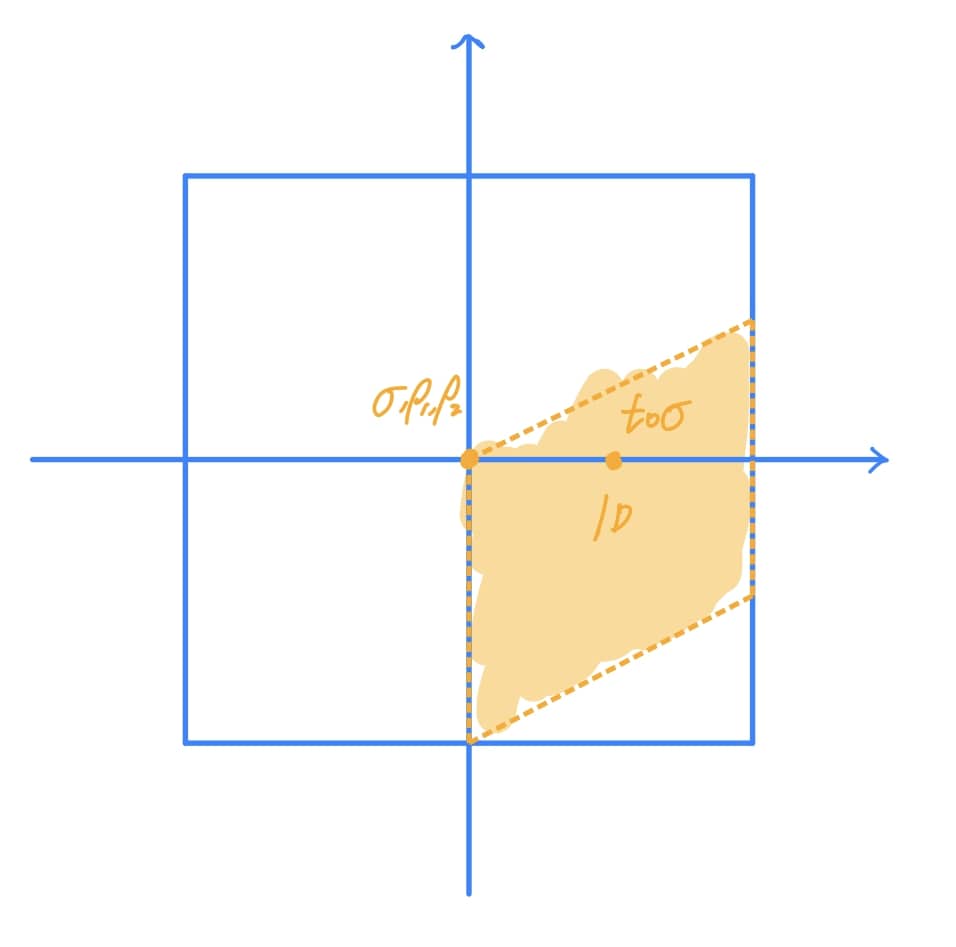}  
      \caption{Fundamental Domain of cmm2}  
      \end{minipage}  
      \end{figure}
  
    We compute $K^{G}_{\bullet}(\underline{E}G)\cong K^{\mathbb{Z}^{2}\rtimes D}_{\bullet}(\mathbb{R}^{2})$.
  
    1. Let $t=(t_{1},0)$ and  $s=(\frac{1}{2}t_{1},\frac{1}{2}s_{2})$ form a basis for $A$ ($t_{1},s_{2}\neq0$ ). 
    
    Point group $D=\{id_{D},\sigma,\rho_{1},\rho_{2}\}$ and the pullback:
    \begin{align*}
      \gamma(\sigma)=((0,0),M_{\sigma}),\quad \gamma(\rho_{1})=((0,0),M_{\rho_{1}}),\quad \gamma(\rho_{2})=((0,0),M_{\rho_{2}}) \\
    M_{\sigma}=\begin{bmatrix}
    -1 & 0\\ 
    0 & -1
    \end{bmatrix}, \quad
    M_{\rho_{1}}=\begin{bmatrix}
      1 & 0\\ 
      1 & -1
      \end{bmatrix}, \quad
      M_{\rho_{2}}=\begin{bmatrix}
        -1 & 0\\ 
        -1 & 1
        \end{bmatrix} .
    \end{align*}
    
    2. $FC(G)=\{[id_{D}],[\sigma],[\rho_{1}],[\rho_{2}],[t\circ \sigma]\}$, where
    $$
\begin{aligned}
t \circ \sigma =((1,0),M_{\sigma})=((\frac{1}{2},0)+(-\frac{1}{2},0)M_{\sigma},M_{\sigma}).
\end{aligned}
$$
    
3. See Table \ref{cmm2}.
 \begin{table}[H]
      \footnotesize
      \centering%
      \begin{tabu} to 0.95\textwidth{X[c]X[c]X[c]X[c]X[c]X[c]}
      \hline
      representative         & identity      & $\sigma$ & $\rho_{1}$ & $\rho_{2}$ &$t\circ\sigma$ \\ \hline
      $X^{g}$                & $\mathbb{R}^{2}$ & $(0,0)$&$(0,0)$&$(0,0)$&$(\frac{1}{2},0)$\\
      $C_{D}(g)$       & $G$    & $D$ &$D$  &$D$ &$t\circ\sigma$  \\
      $X^{g}/C_{D}(g)$ & Fundamental Domain & $(0,0)$&$(0,0)$&$(0,0)$&$(\frac{1}{2},0)$\\
      \tiny{$\bigoplus \limits_{k \in even} H_{k}(X^{g}/C_{D}(g))$} & $\mathbb{Z}$ & $\mathbb{Z}$ & $\mathbb{Z}$ &$\mathbb{Z}$&$\mathbb{Z}$ \\
      \tiny{$\bigoplus \limits_{k \in odd} H_{k}(X^{g}/C_{D}(g))$} &
      0&0&0&0&0 \\ \hline
      \end{tabu}%
      \caption{The computation about cmm2}
      \label{cmm2}
      \end{table}
    
      4. $K_{0}(cmm2)=\mathbb{Z}^{5}, K_{1}(cmm2)=0$.

\subsection*{pmm2}
\begin{figure}[btbp]
  \centering
  \begin{minipage}[t]{0.45\textwidth}   
  \centering
  \includegraphics[height=5cm]{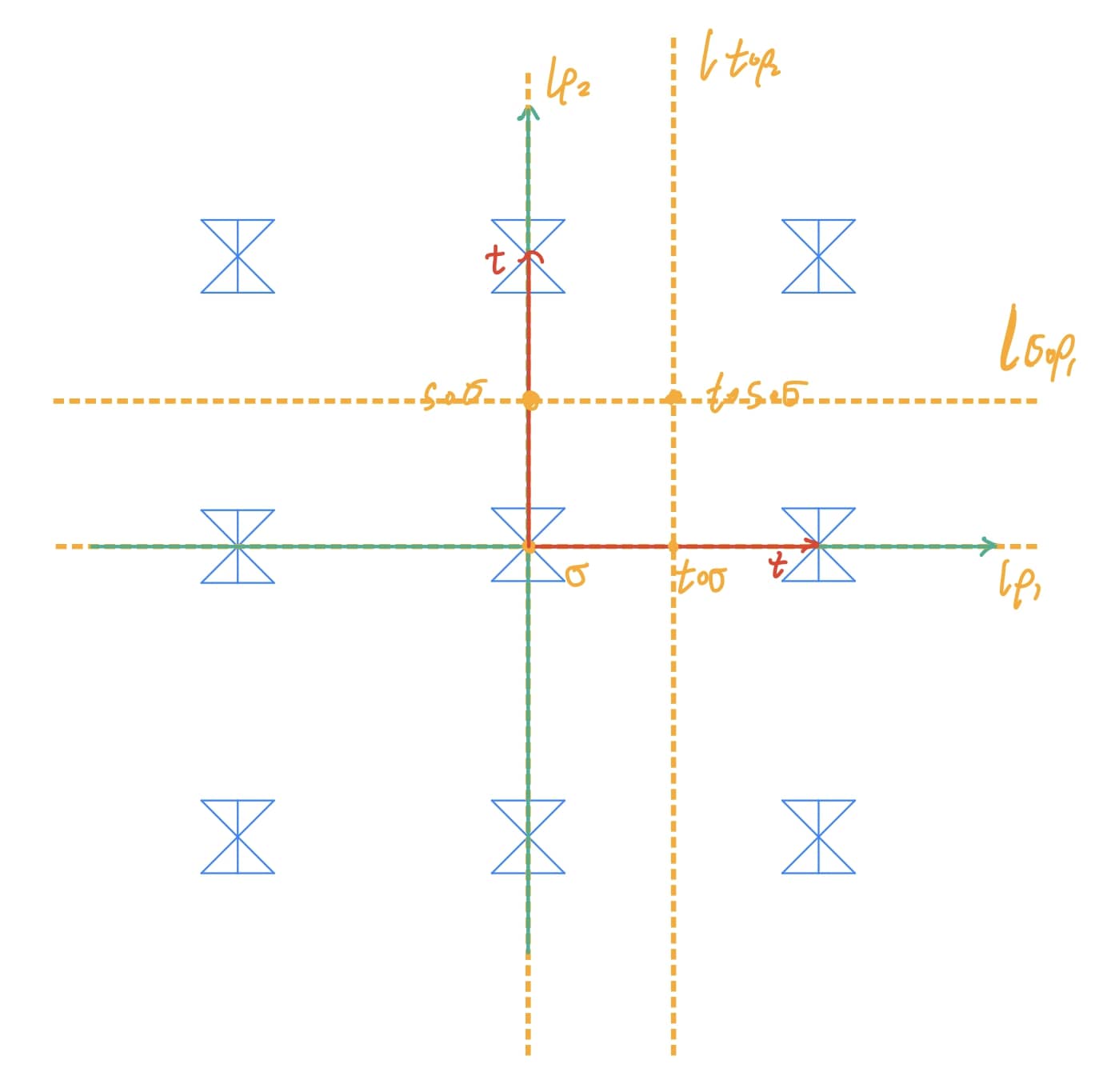}
  \caption{pmm2}  
  \end{minipage}
  \begin{minipage}[t]{0.45\textwidth}  
  \centering  
  \includegraphics[height=5cm]{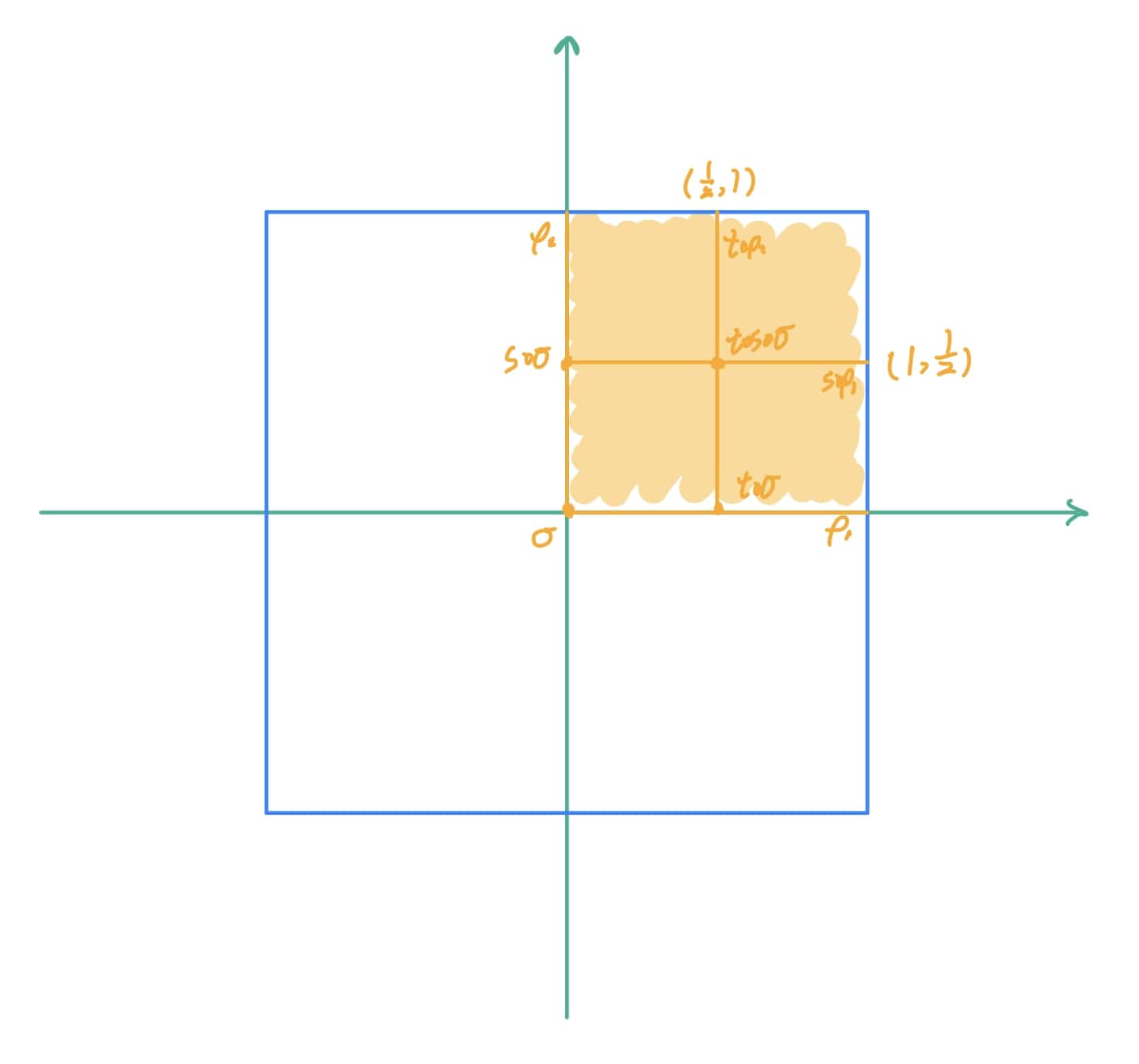}  
  \caption{Fundamental Domain of pmm2}  
  \end{minipage}  
  \end{figure}

We compute $K^{G}_{\bullet}(\underline{E}G)\cong K^{\mathbb{Z}^{2}\rtimes D}_{\bullet}(\mathbb{R}^{2})$.

1. Let $t=(t_{1},0)$ and  $s=(0,s_{2})$ form a basis for $A$ ($t_{1},s_{2}\neq0$ ). 

Point group $D=\{id_{D},\sigma,\rho_{1},\rho_{2}\}$ and the pullback:
\begin{align*}
  \gamma(\sigma)=((0,0),M_{\sigma}),\quad \gamma(\rho_{1})=((0,0),M_{\rho_{1}}),\quad \gamma(\rho_{2})=((0,0),M_{\rho_{2}}) \\
M_{\sigma}=\begin{bmatrix}
-1 & 0\\ 
0 & -1
\end{bmatrix}, \quad
M_{\rho_{1}}=\begin{bmatrix}
  1 & 0\\ 
  0 & -1
  \end{bmatrix}, \quad
  M_{\rho_{2}}=\begin{bmatrix}
    -1 & 0\\ 
    0 & 1
    \end{bmatrix} .
\end{align*}

2. $FC(G)=\{[id_{D}],[\sigma],[\rho_{1}],[\rho_{2}],[t\circ \sigma],[s\circ\sigma],[t\circ s\circ\sigma],[t\circ\rho_{2}],[s\circ\rho_{1}]\}$, where
$$
\begin{aligned}
t \circ \sigma =((1,0),M_{\sigma})=((\frac{1}{2},0)+(-\frac{1}{2},0)M_{\sigma},M_{\sigma}), \\
s \circ \sigma =((0,1),M_{\sigma})=((0,\frac{1}{2})+(0,-\frac{1}{2})M_{\sigma},M_{\sigma}),\\
s\circ t \circ \sigma =((1,1),M_{\sigma})=((\frac{1}{2},\frac{1}{2})+(-\frac{1}{2},-\frac{1}{2})M_{\sigma},M_{\sigma}),\\
t \circ \rho_{2} =((1,0),M_{\rho_{2}})=((\frac{1}{2},0)+(-\frac{1}{2},0)M_{\rho_{2}},M_{\rho_{2}}),\\
s \circ \rho_{1} =((0,1),M_{\rho_{1}})=((0,\frac{1}{2})+(0,-\frac{1}{2})M_{\rho_{1}},M_{\rho_{1}}).
\end{aligned}
$$

3. See Table \ref{pmm2}.
\begin{table}[H]
  \footnotesize
  \centering%
  \begin{tabu} to 0.95\textwidth{X[c]X[c]X[c]X[c]X[c]}
  \hline
  representative         & identity      & $\rho_{1}$ & $\rho_{2}$ & $\sigma$  \\ \hline
  $X^{g}$                & $\mathbb{R}^{2}$ & $(x,0)$&$(0,y)$&$(0,0)$\\
  $C_{D}(g)$       & $G$    & $t,\sigma,\rho_{1},\rho_{2}$ &$s,\sigma,\rho_{1},\rho_{2}$  &$D$   \\
  $X^{g}/C_{D}(g)$ & Fundamental Domain & $\{(x,0)|0\leqslant x\leqslant 1\}$&$\{(0,y)|0\leqslant y\leqslant 1\}$&$(0,0)$ \\
  \tiny{$\bigoplus \limits_{k \in even} H_{k}(X^{g}/C_{D}(g))$} & $\mathbb{Z}$ & $\mathbb{Z}$ & $\mathbb{Z}$ &$\mathbb{Z}$\\
  \tiny{$\bigoplus \limits_{k \in odd} H_{k}(X^{g}/C_{D}(g))$} &
  0&0&0&0 \\ \hline \hline
  $t\circ\sigma$ & $s\circ\sigma$&$s\circ t\circ\sigma$&$t\circ\rho_{2}$&$s\circ\rho_{1}$ \\ \hline
  $(0,0)$&$(0,\frac{1}{2})$&$(\frac{1}{2},\frac{1}{2})$&$(\frac{1}{2},y)$&$(x,\frac{1}{2})$ \\
  $T\circ\sigma$&$\rho_{1},\rho_{2},\sigma$&$\rho_{1},\rho_{2},\sigma$&$t,\rho_{1},\rho_{2},\sigma$&$s,\rho_{1},\rho_{2},\sigma$ \\
  $(0,0)$&$(0,\frac{1}{2})$&$(\frac{1}{2},\frac{1}{2})$&$\{(\frac{1}{2},y)|0\leqslant y\leqslant 1\}$&$\{(x,\frac{1}{2})|0\leqslant x\leqslant 1\}$ \\ 
  $\mathbb{Z}$ & $\mathbb{Z}$ & $\mathbb{Z}$ &$\mathbb{Z}$ &$\mathbb{Z}$\\
  0&0&0&0&0 \\
  \hline
  \end{tabu}%
  \caption{The computation about pmm2}
  \label{pmm2}
  \end{table}

  4. $K_{0}(pmm2)=\mathbb{Z}^{9}, K_{1}(pmm2)=0$.
  \subsection*{pmg}
\begin{figure}[btbp]
  \centering
  \begin{minipage}[t]{0.45\textwidth}   
  \centering
  \includegraphics[height=5cm]{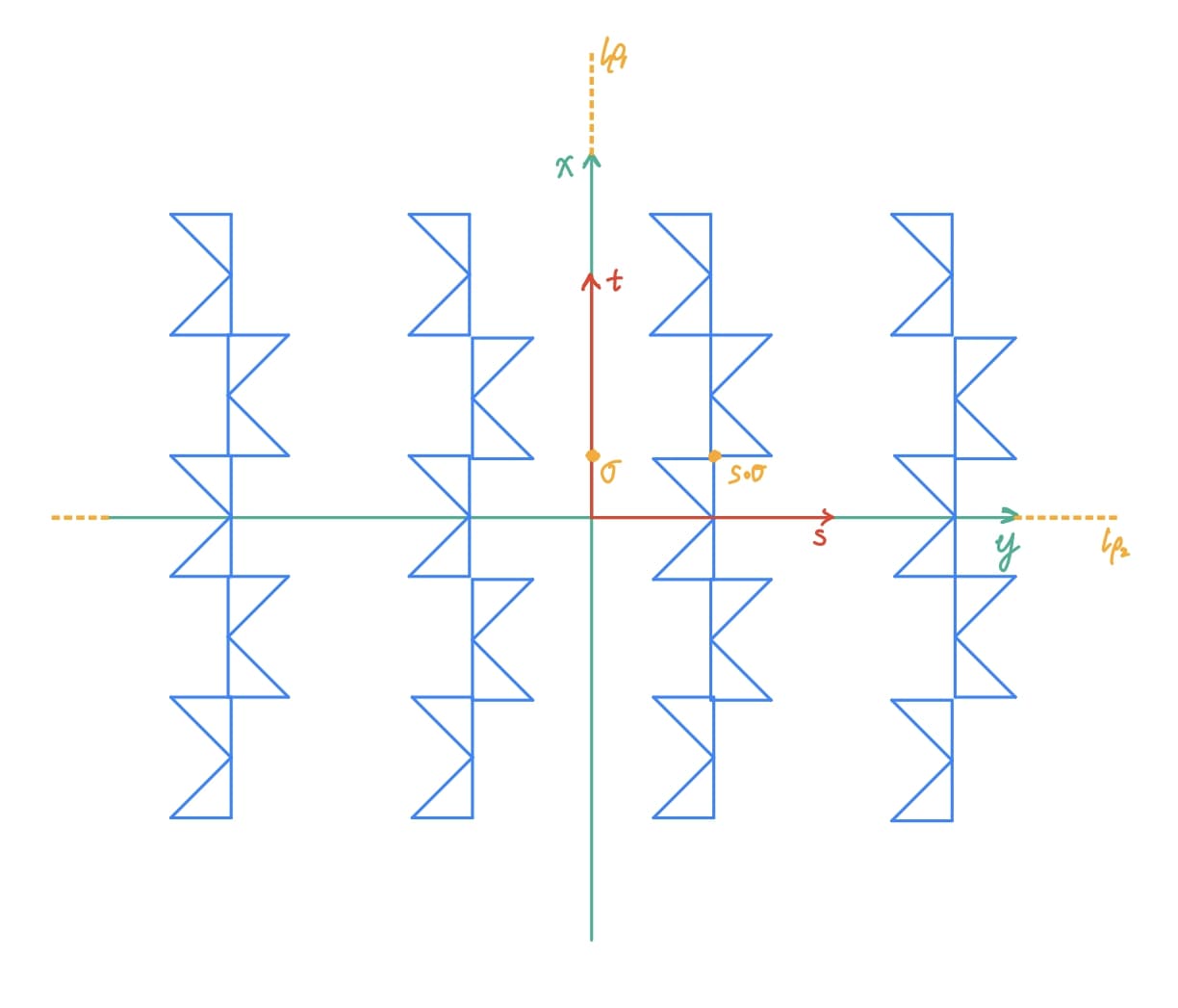}
  \caption{pmg}  
  \end{minipage}
  \begin{minipage}[t]{0.45\textwidth}  
  \centering  
  \includegraphics[height=5cm]{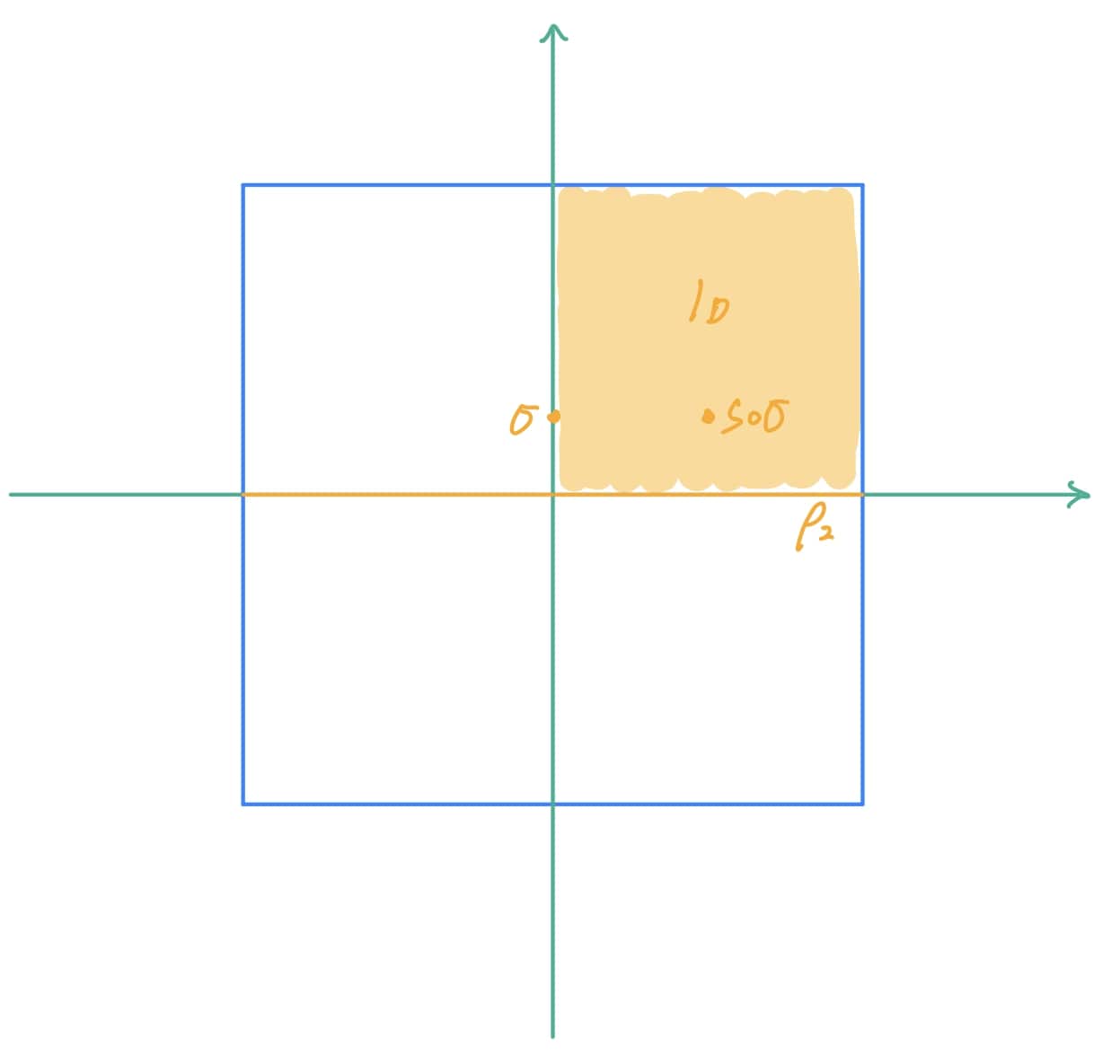}  
  \caption{Fundamental Domain of pmg}  
  \end{minipage}  
  \end{figure}

We compute $K^{G}_{\bullet}(\underline{E}G)\cong K^{\mathbb{Z}^{2}\rtimes D}_{\bullet}(\mathbb{R}^{2})$.

1. Let $t=(t_{1},0)$ and  $s=(0,s_{2})$ form a basis for $A$ ($t_{1},s_{2}\neq0$ ). 

Point group $D=\{id_{D},\sigma,\rho_{1},\rho_{2}\}$ and the pullback:
\begin{align*}
  \gamma(\sigma)=((\frac{1}{2},0),M_{\sigma}),\quad \gamma(\rho_{1})=((\frac{1}{2},0),M_{\rho_{1}}),\quad \gamma(\rho_{2})=((0,0),M_{\rho_{2}}) \\
M_{\sigma}=\begin{bmatrix}
-1 & 0\\ 
0 & -1
\end{bmatrix}, \quad
M_{\rho_{1}}=\begin{bmatrix}
  1 & 0\\ 
  0 & -1
  \end{bmatrix}, \quad
  M_{\rho_{2}}=\begin{bmatrix}
    -1 & 0\\ 
    0 & 1
    \end{bmatrix} .
\end{align*}

2. $FC(G)=\{[id_{D}],[\sigma],[\rho_{2}],[s\circ\sigma]$, where
$$
\begin{aligned}
s \circ \sigma =((\frac{1}{2},1),M_{\sigma})=((\frac{1}{4},\frac{1}{2})+(-\frac{1}{4},-\frac{1}{2})M_{\sigma},M_{\sigma}).
\end{aligned}
$$

3. See Table \ref{pmg}.
\begin{table}[H]
  \footnotesize
  \centering%
  \begin{tabu} to 0.95\textwidth{X[c]X[c]X[c]X[c]X[c]}
  \hline
  representative         & identity      & $\rho_{2}$ & $\sigma$ & $s\circ\sigma$  \\ \hline
  $X^{g}$                & $\mathbb{R}^{2}$ & $(x,0)$&$(\frac{1}{4},0)$&$(\frac{1}{4},\frac{1}{2})$\\
  $C_{D}(g)$       & $G$    & $t,id_{D},\rho_{2}$ &$\sigma,id_{D}$  &  $id_{D},s\circ\sigma$  \\
  $X^{g}/C_{D}(g)$ & Fundamental Domain & $\{(x,0)|-1\leqslant x\leqslant 1\}$&$(\frac{1}{4},0)$&$(\frac{1}{4},\frac{1}{2})$ \\
  \tiny{$\bigoplus \limits_{k \in even} H_{k}(X^{g}/C_{D}(g))$} & $\mathbb{Z}$ & $\mathbb{Z}$ & $\mathbb{Z}$ &$\mathbb{Z}$\\
  \tiny{$\bigoplus \limits_{k \in odd} H_{k}(X^{g}/C_{D}(g))$} &
  0&$\mathbb{Z}$&0&0 \\ 
  \hline
  \end{tabu}%
  \caption{The computation about pmg}
  \label{pmg}
  \end{table}

  4. $K_{0}(pmg)=\mathbb{Z}^{4}, K_{1}(pmg)=\mathbb{Z}$.

\subsection*{pgg2}
\begin{figure}[btbp]
  \centering
  \begin{minipage}[t]{0.45\textwidth}   
  \centering
  \includegraphics[height=5cm]{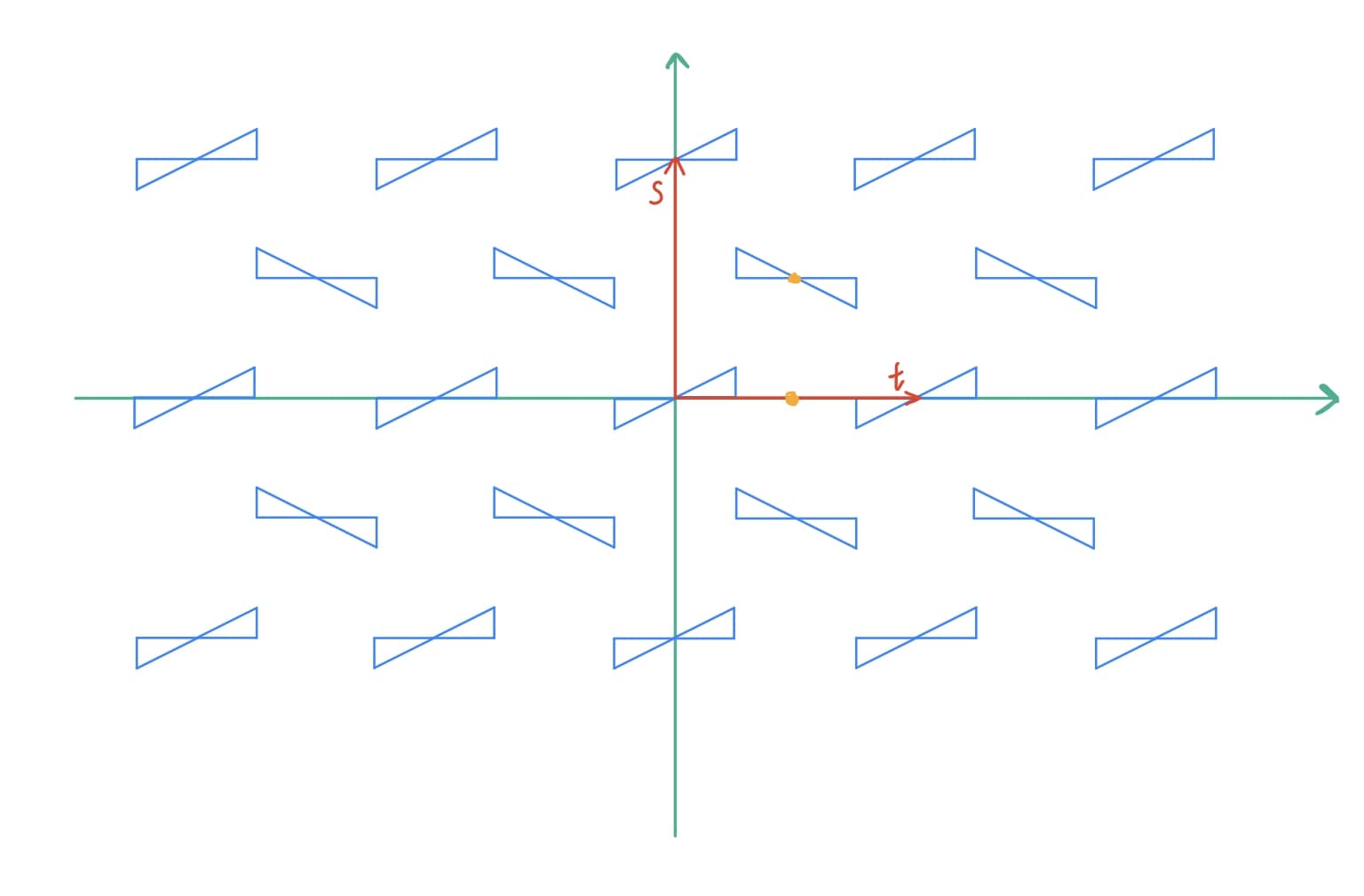}
  \caption{pgg2}  
  \end{minipage}
  \begin{minipage}[t]{0.45\textwidth}  
  \centering  
  \includegraphics[height=5cm]{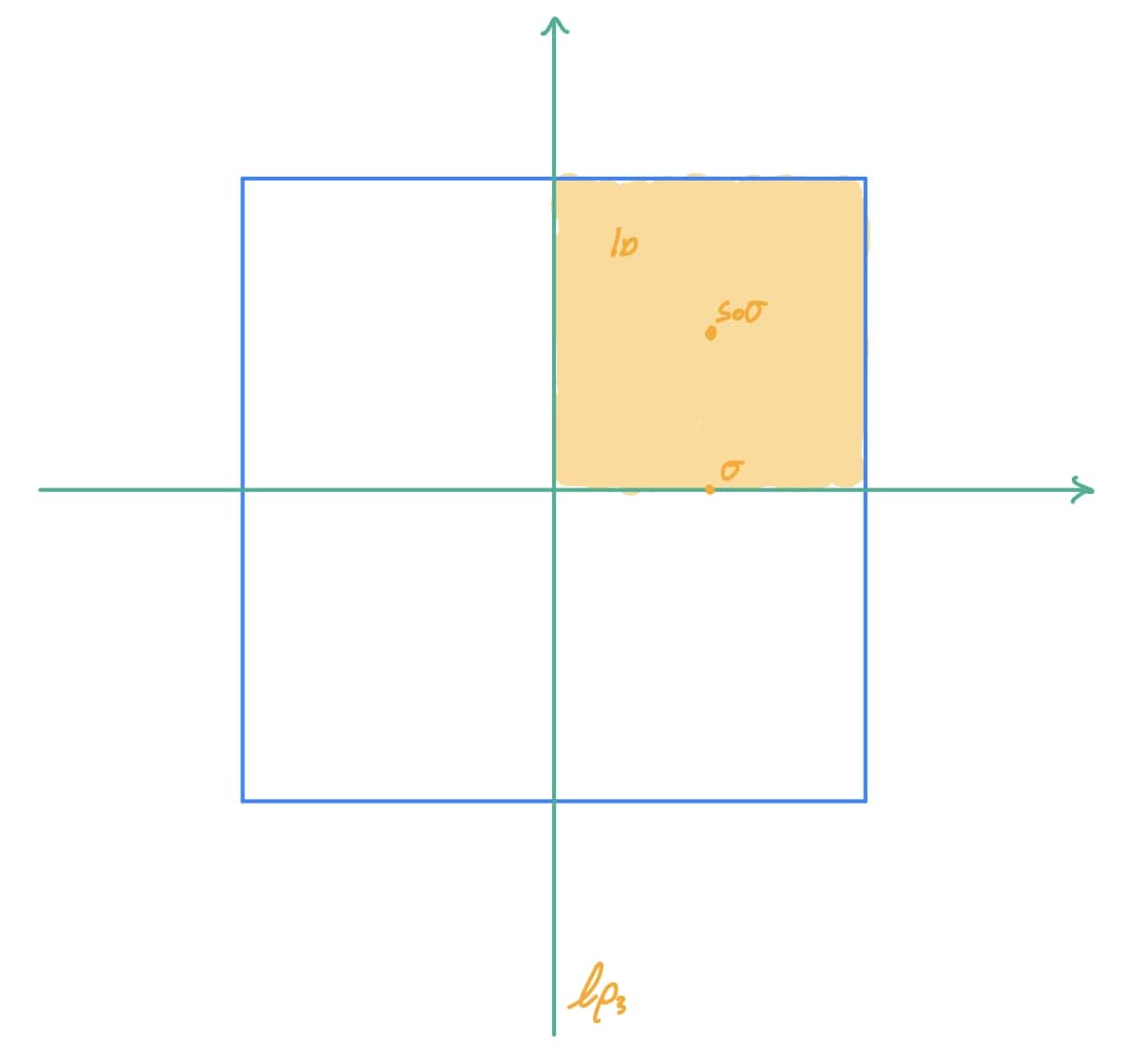}  
  \caption{Fundamental Domain of pgg2}  
  \end{minipage}  
  \end{figure}

We compute $K^{G}_{\bullet}(\underline{E}G)\cong K^{\mathbb{Z}^{2}\rtimes D}_{\bullet}(\mathbb{R}^{2})$.

1. Let $t=(t_{1},0)$ and  $s=(0,s_{2})$ form a basis for $A$ ($t_{1},s_{2}\neq0$ ). 

Point group $D=\{id_{D},\sigma,\rho_{1},\rho_{2}\}$ and the pullback:
\begin{align*}
  \gamma(\sigma)=((1,0),M_{\sigma}),\quad \gamma(\rho_{1})=((\frac{1}{2},\frac{1}{2}),M_{\rho_{1}}),\quad \gamma(\rho_{2})=((\frac{1}{2},\frac{1}{2}),M_{\rho_{2}}) \\
M_{\sigma}=\begin{bmatrix}
-1 & 0\\ 
0 & -1
\end{bmatrix}, \quad
M_{\rho_{1}}=\begin{bmatrix}
  1 & 0\\ 
  0 & -1
  \end{bmatrix}, \quad
  M_{\rho_{2}}=\begin{bmatrix}
    -1 & 0\\ 
    0 & 1
    \end{bmatrix} .
\end{align*}

2. $FC(G)=\{[id_{D}],[\sigma],[s\circ\sigma]$, where
$$
\begin{aligned}
s \circ \sigma =((1,1),M_{\sigma})=((\frac{1}{2},\frac{1}{2})+(-\frac{1}{2},-\frac{1}{2})M_{\sigma},M_{\sigma}).
\end{aligned}
$$

3. See Table \ref{pgg2}.
\begin{table}[H]
  \footnotesize
  \centering%
  \begin{tabu} to 0.95\textwidth{X[c]X[c]X[c]X[c]}
  \hline
  representative         & identity       & $\sigma$ & $s\circ\sigma$  \\ \hline
  $X^{g}$                & $\mathbb{R}^{2}$ &$(\frac{1}{2},0)$&$(\frac{1}{2},\frac{1}{2})$\\
  $C_{D}(g)$       & $G$     &$\sigma,id_{D}$  &  $id_{D},s\circ\sigma$  \\
  $X^{g}/C_{D}(g)$ & Fundamental Domain &$(\frac{1}{2},0)$&$(\frac{1}{2},\frac{1}{2})$ \\
  \tiny{$\bigoplus \limits_{k \in even} H_{k}(X^{g}/C_{D}(g))$} & $\mathbb{Z}$ & $\mathbb{Z}$  &$\mathbb{Z}$\\
  \tiny{$\bigoplus \limits_{k \in odd} H_{k}(X^{g}/C_{D}(g))$} &
  $\mathbb{Z}_{2}$&0&0 \\ 
  \hline
  \end{tabu}%
  \caption{The computation about pgg2}
  \label{pgg2}
  \end{table}

  4. $K_{0}(pgg2)=\mathbb{Z}^{3}, K_{1}(pgg2)=\mathbb{Z}_{2}$.

  \subsection*{p31m}
    \begin{figure}[btbp]
      \centering
      \begin{minipage}[t]{0.45\textwidth}   
      \centering
      \includegraphics[height=5cm]{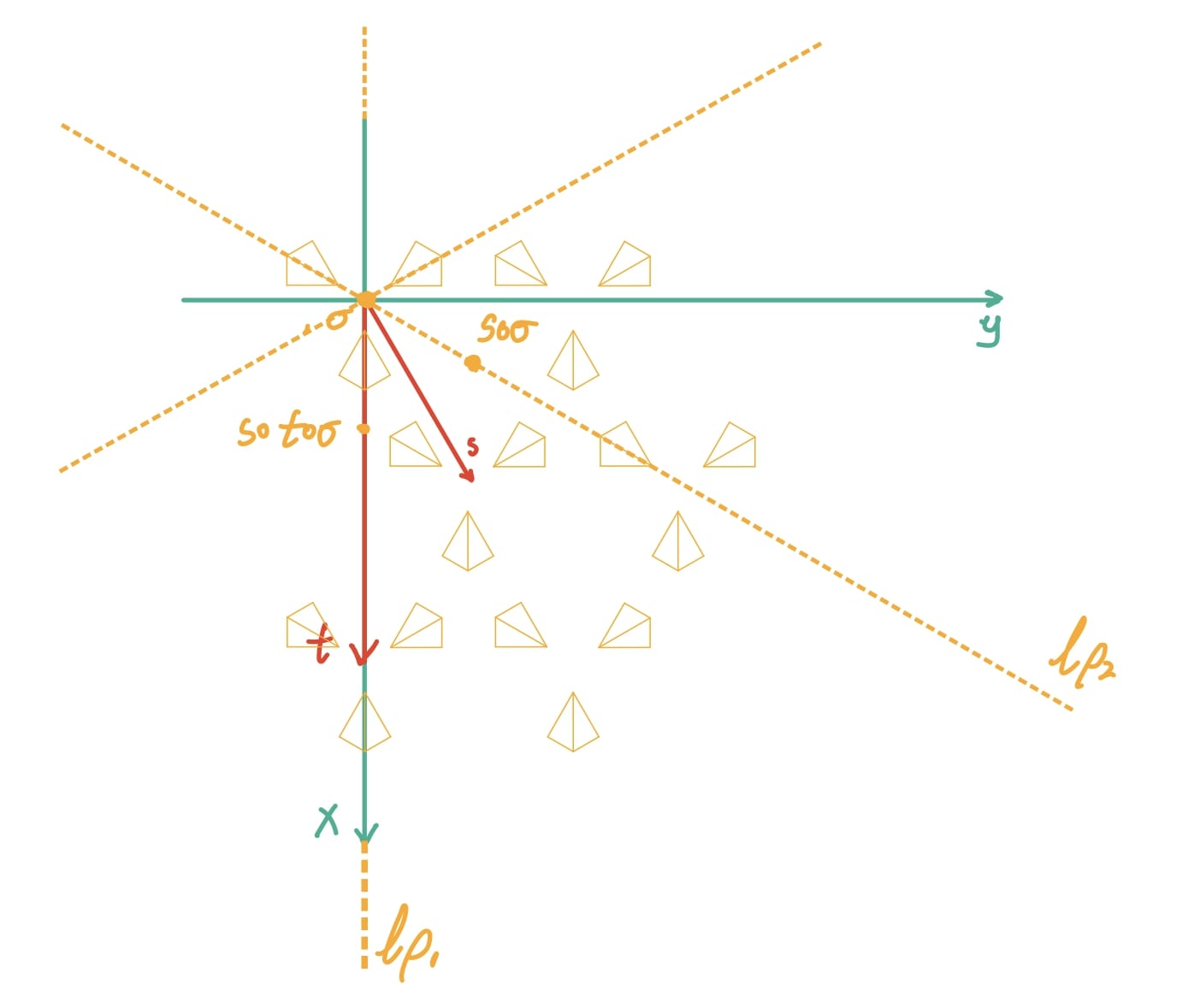}
      \caption{p31m}  
      \end{minipage}
      \begin{minipage}[t]{0.45\textwidth}  
      \centering  
      \includegraphics[height=5cm]{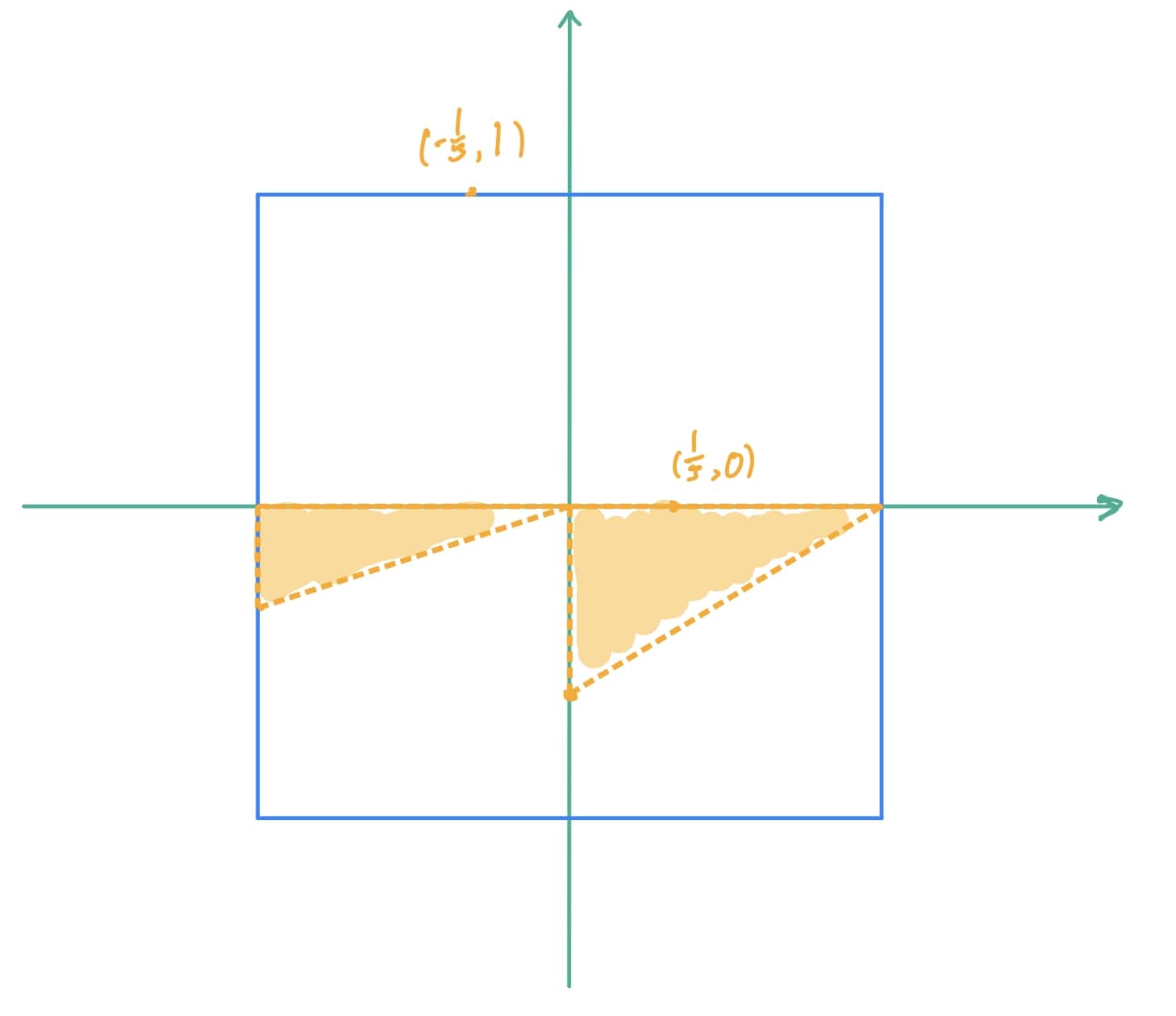}  
      \caption{Fundamental Domain of p31m}  
      \end{minipage}  
      \end{figure}
  
    We compute $K^{G}_{\bullet}(\underline{E}G)\cong K^{\mathbb{Z}^{2}\rtimes D}_{\bullet}(\mathbb{R}^{2})$.
  
    1. Let $t=(t_{1},0)$ and  $s=(\frac{1}{2}t_{1},\frac{\sqrt{3}}{6}s_{2})$ form a basis for $A$ ($t_{1},s_{2}\neq0$ ). 
    
    Point group $D=\{id_{D},\sigma,\sigma^{2},\rho_{1},\rho_{2},\rho_{3}\}$ and the pullback of the genetors:
    \begin{align*}
      \gamma(\sigma)=((0,0),M_{\sigma}),\quad \gamma(\rho_{1})=((0,0),M_{\rho_{1}}), \\
    M_{\sigma}=\begin{bmatrix}
    -2 & 3\\ 
    -1 & 1
    \end{bmatrix}, \quad
    M_{\rho_{1}}=\begin{bmatrix}
      1 & 0\\ 
      1 & -1
      \end{bmatrix} .
    \end{align*}
    
    2. $FC(G)=\{[id_{D}],[\sigma],[\rho_{1}],[-s\circ\sigma],[-s\circ t\circ \sigma]\}$, where
    $$
\begin{aligned}
s \circ \sigma =((0,1),M_{\sigma})=((-\frac{1}{3},1)+(\frac{1}{3},-1)M_{\sigma},M_{\sigma}),\\
s \circ t\circ \sigma =((1,1),M_{\sigma})=((\frac{1}{3},0)+(-\frac{1}{3},0)M_{\sigma},M_{\sigma}).
\end{aligned}
$$
    
3. See Table \ref{p31m}.
 \begin{table}[H]
      \footnotesize
      \centering%
      \begin{tabu} to 0.95\textwidth{X[c]X[c]X[c]X[c]X[c]X[c]}
      \hline
      representative         & identity      & $\rho_{1}$ & $\sigma$ & $s\circ\sigma$ &$s\circ t\circ\sigma$ \\ \hline
      $X^{g}$                & $\mathbb{R}^{2}$ & $(x,0)$&$(0,0)$&$(-\frac{1}{3},1)$&$(\frac{1}{3},0)$\\
      $C_{D}(g)$       & $G$    & $t,\rho_{1}$ &$id_{D},\sigma,\sigma^{2}$  &$id_{D},s\circ\sigma$ &$id_{D},s\circ t\circ\sigma$  \\
      $X^{g}/C_{D}(g)$ & Fundamental Domain & $\{(x,0)|-1\leqslant x\leqslant 1\}$&$(0,0)$&$(-\frac{1}{3},0),0)$&$(\frac{1}{3},0)$\\
      \tiny{$\bigoplus \limits_{k \in even} H_{k}(X^{g}/C_{D}(g))$} & $\mathbb{Z}$ & $\mathbb{Z}$ & $\mathbb{Z}$ &$\mathbb{Z}$&$\mathbb{Z}$ \\
      \tiny{$\bigoplus \limits_{k \in odd} H_{k}(X^{g}/C_{D}(g))$} &
      0&$\mathbb{Z}$&0&0&0 \\ \hline
      \end{tabu}%
      \caption{The computation about p31m}
      \label{p31m}
      \end{table}
    
  4. $K_{0}(p31m)=\mathbb{Z}^{5}, K_{1}(p31m)=\mathbb{Z}$.

  \subsection*{p3m1}
  \begin{figure}[btbp]
    \centering
    \begin{minipage}[t]{0.45\textwidth}   
    \centering
    \includegraphics[height=5cm]{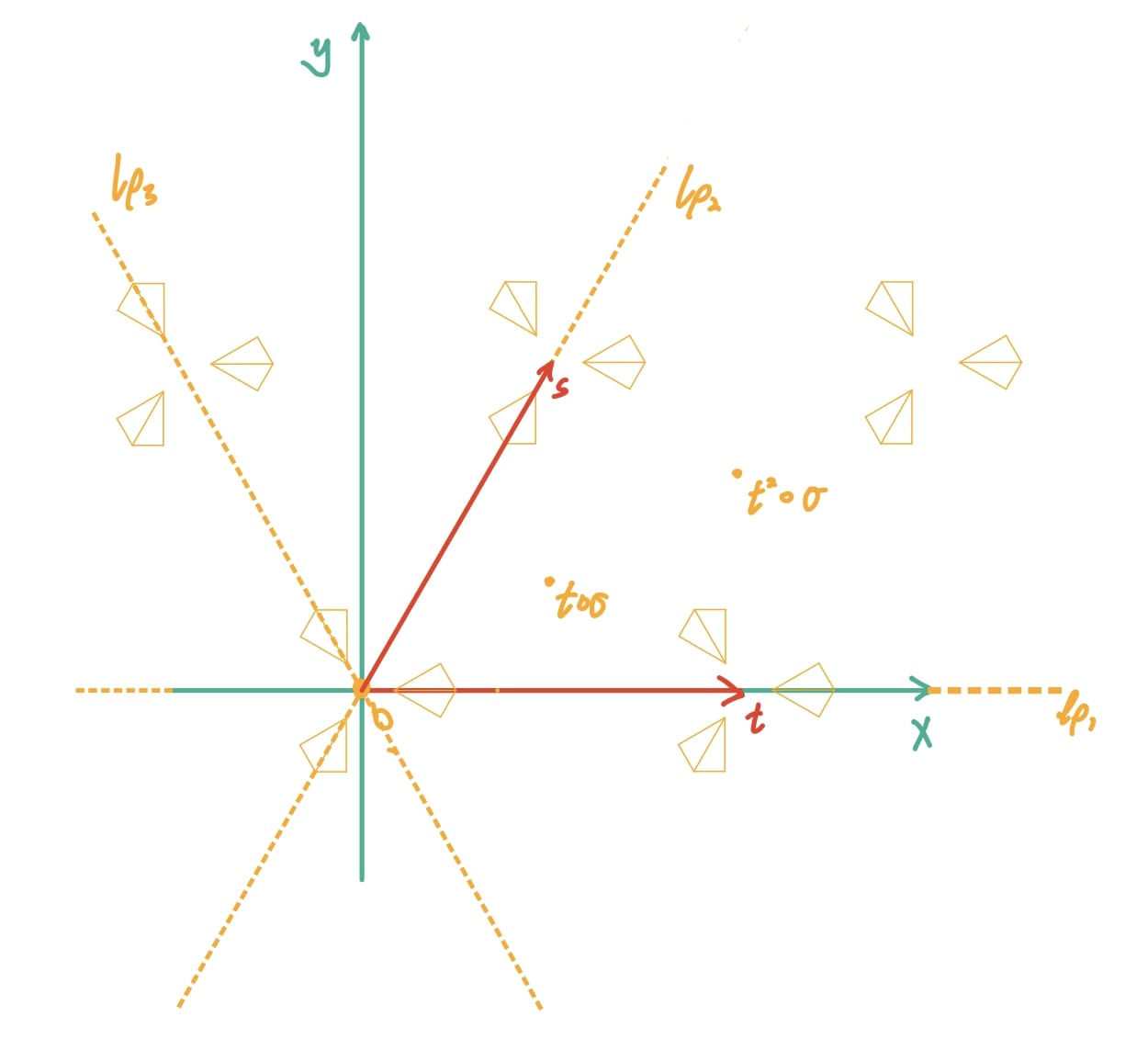}
    \caption{p3m1}  
    \end{minipage}
    \begin{minipage}[t]{0.45\textwidth}  
    \centering  
    \includegraphics[height=5cm]{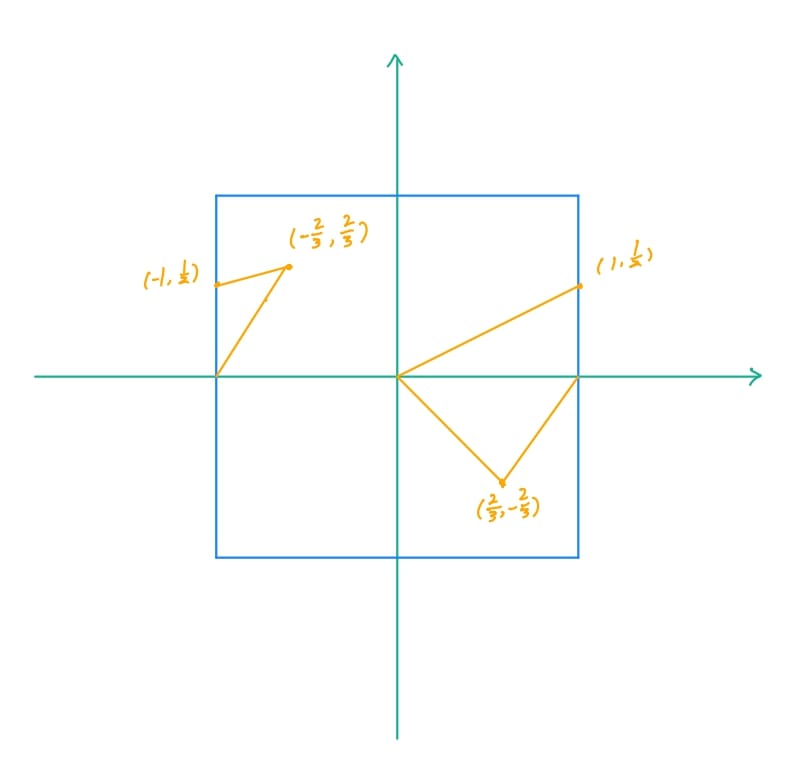}  
    \caption{Fundamental Domain of p3m1}  
    \end{minipage}  
    \end{figure}

  We compute $K^{G}_{\bullet}(\underline{E}G)\cong K^{\mathbb{Z}^{2}\rtimes D}_{\bullet}(\mathbb{R}^{2})$.

  1. Let $t=(t_{1},0)$ and  $s=(\frac{1}{2}t_{1},\frac{\sqrt{3}}{2}s_{2})$ form a basis for $A$ ($t_{1},s_{2}\neq0$ ). 
  
  Point group $D=\{id_{D},\sigma,\sigma^{2},\rho_{1},\rho_{2},\rho_{3}\}$ and the pullback of the genetors:
  \begin{align*}
    \gamma(\sigma)=((0,0),M_{\sigma}),\quad \gamma(\rho_{1})=((0,0),M_{\rho_{1}}), \\
  M_{\sigma}=\begin{bmatrix}
  -1 & 1\\ 
  -1 & 0
  \end{bmatrix}, \quad
  M_{\rho_{1}}=\begin{bmatrix}
    1 & 0\\ 
    1 & -1
    \end{bmatrix} .
  \end{align*}
  
  2. $FC(G)=\{[id_{D}],[\sigma],[\rho_{1}],[t\circ\sigma],[2t\circ t\circ \sigma]\}$, where
  $$
\begin{aligned}
t \circ \sigma =((1,0),M_{\sigma})=((\frac{1}{3},\frac{1}{3})+(-\frac{1}{3},-\frac{1}{3})M_{\sigma},M_{\sigma}),\\
2t \circ t\circ \sigma =((2,0),M_{\sigma})=((\frac{2}{3},\frac{2}{3})+(-\frac{2}{3},-\frac{2}{3})M_{\sigma},M_{\sigma}).
\end{aligned}
$$
  
3. See Table \ref{p3m1}.
\begin{table}[H]
    \footnotesize
    \centering%
    \begin{tabu} to 0.95\textwidth{X[c]X[c]X[c]X[c]X[c]X[c]}
    \hline
    representative         & identity      & $\rho_{1}$ & $\sigma$ & $t\circ\sigma$ &$2t\circ\sigma$ \\ \hline
    $X^{g}$                & $\mathbb{R}^{2}$ & $(0,y)$&$(0,0)$&$(\frac{1}{3},\frac{1}{3})$&$(\frac{2}{3},\frac{2}{3})$\\
    $C_{D}(g)$       & $G$    & $t,\rho_{1}$ &$id_{D},\sigma,\sigma^{2}$  &$id_{D},t\circ\sigma$ &$id_{D}, 2t\circ\sigma$  \\
    $X^{g}/C_{D}(g)$ & Fundamental Domain & $\{(0,y)|-1\leqslant y\leqslant 1\}$&$(0,0)$&$(\frac{1}{3},\frac{1}{3})$&$(\frac{2}{3},\frac{2}{3})$\\
    \tiny{$\bigoplus \limits_{k \in even} H_{k}(X^{g}/C_{D}(g))$} & $\mathbb{Z}$ & $\mathbb{Z}$ & $\mathbb{Z}$ &$\mathbb{Z}$&$\mathbb{Z}$ \\
    \tiny{$\bigoplus \limits_{k \in odd} H_{k}(X^{g}/C_{D}(g))$} &
    0&$\mathbb{Z}$&0&0&0 \\ \hline
    \end{tabu}%
    \caption{The computation about p3m1}
    \label{p3m1}
    \end{table}
  
4. $K_{0}(p3m1)=\mathbb{Z}^{5}, K_{1}(p3m1)=\mathbb{Z}$.

\subsection*{p4mm}
\begin{figure}[btbp]
  \centering
  \begin{minipage}[t]{0.45\textwidth}   
  \centering
  \includegraphics[height=5cm]{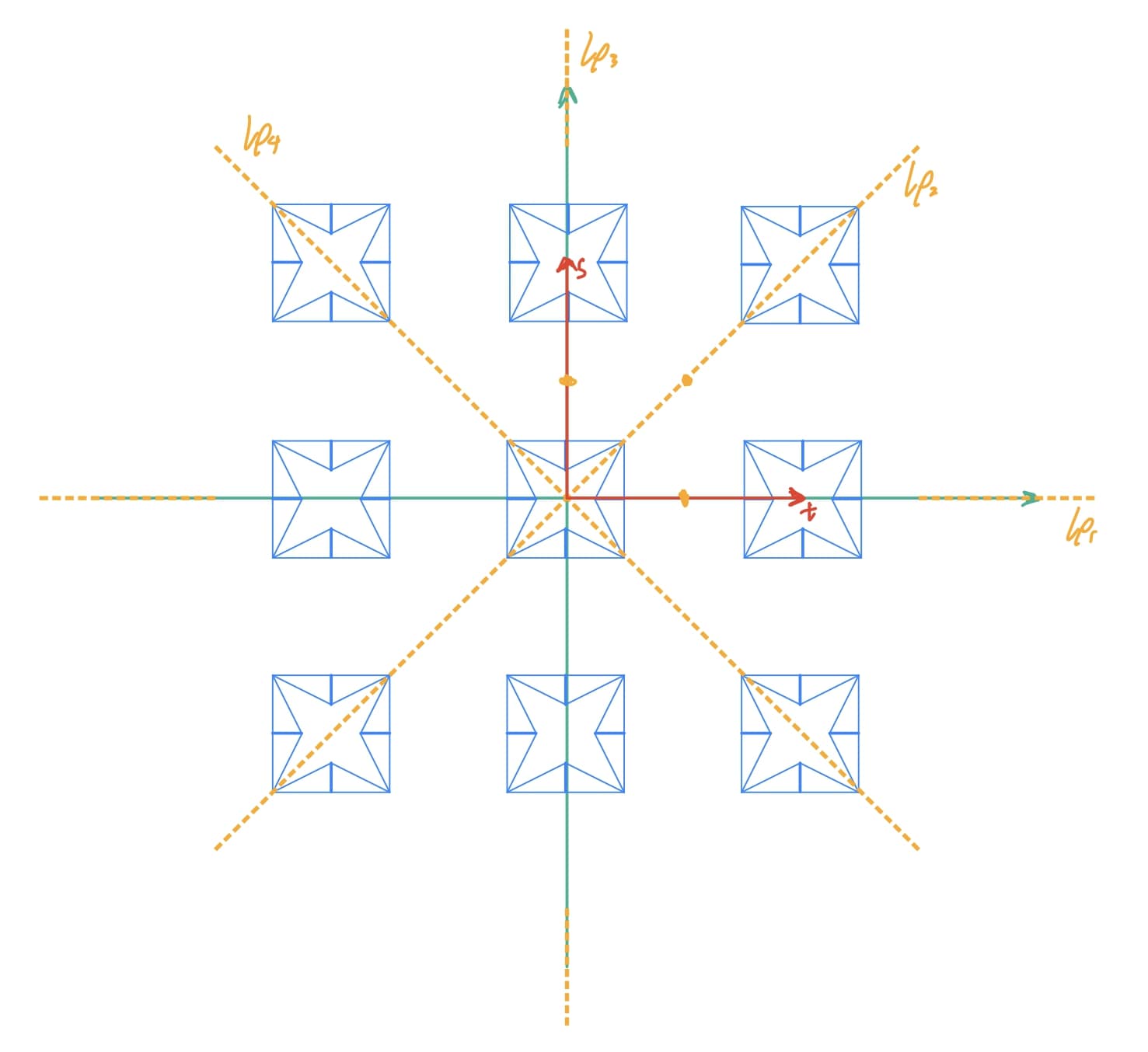}
  \caption{p4mm}  
  \end{minipage}
  \begin{minipage}[t]{0.45\textwidth}  
  \centering  
  \includegraphics[height=5cm]{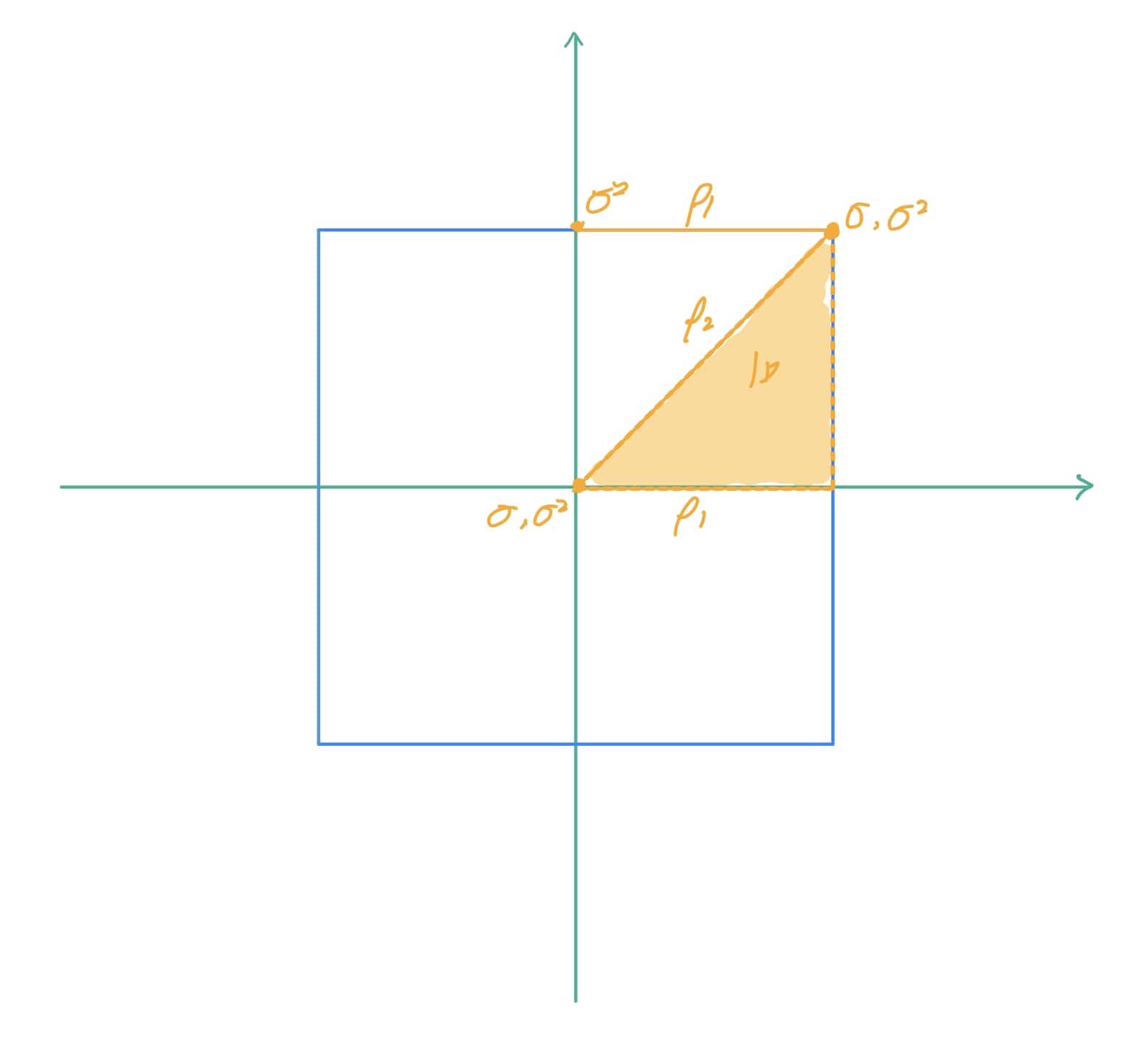}  
  \caption{Fundamental Domain of p4mm}  
  \end{minipage}  
  \end{figure}

We compute $K^{G}_{\bullet}(\underline{E}G)\cong K^{D}_{\bullet}(\mathbb{T}^{2})$.

1. Let $t=(t_{1},0)$ and  $s=(0,s_{2})$ form a basis for $A$ ($t_{1},s_{2}\neq0$ ). 

Point group $D=\{id_{D},\sigma,\sigma^{2},\sigma^{3},\rho_{1},\rho_{2},\rho_{3},\rho_{4}\}$ and the pullback of the genetors:
\begin{align*}
  \gamma(\sigma)=((0,0),M_{\sigma}),\quad \gamma(\rho_{1})=((0,0),M_{\rho_{1}}), \\
M_{\sigma}=\begin{bmatrix}
0 & 1\\ 
-1 & 0
\end{bmatrix}, \quad
M_{\rho_{1}}=\begin{bmatrix}
  1 & 0\\ 
  0 & -1
  \end{bmatrix} .
\end{align*}

2. $FC(G)=\{[id_{D}],[\sigma],[\sigma^{2}],[\rho_{1}],[\rho_{2}]\}$.

3. See Table \ref{p4mm}.
\begin{table}[H]
  \footnotesize
  \centering%
  \begin{tabu} to 0.95\textwidth{X[c]X[c]X[c]X[c]X[c]X[c]}
  \hline
  representative         & identity      & $\sigma$ & $\sigma^{2}$ & $\rho_{1}$ &$\rho_{2}$ \\ \hline
  $X^{g}$                & $\mathbb{T}^{2}$ & $(\pm 1,\pm 1)$, $(0,0)$ &$(\pm 1,\pm 1)$, $(0,0)$, $(0,\pm 1)$, $(\pm 1,0)$&$(x,0)$, $(x,1)$&$(x,x)$\\
  $C_{D}(g)$       & $D$    & $D$ &$D$  &$D$ &$D$  \\
  $X^{g}/C_{D}(g)$ & Fundamental Domain & $(0,0),(1,1)$&$(0,0)$, $(1,1)$, $(0,1)$&$\{(x,0)|0\leqslant x\leqslant 1\}$, $\{(x,1)|0\leqslant x\leqslant 1\}$&$\{(x,x)|0\leqslant x\leqslant 1\}$\\
  \tiny{$\bigoplus \limits_{k \in even} H_{k}(X^{g}/C_{D}(g))$} & $\mathbb{Z}$ & $\mathbb{Z}^{2}$ & $\mathbb{Z}^{3}$ &$\mathbb{Z}^{2}$&$\mathbb{Z}$ \\
  \tiny{$\bigoplus \limits_{k \in odd} H_{k}(X^{g}/C_{D}(g))$} &
  0&0&0&0&0 \\ \hline
  \end{tabu}%
  \caption{The computation about p4mm}
  \label{p4mm}
  \end{table}

4. $K_{0}(p4mm)=\mathbb{Z}^{9}, K_{1}(p4mm)=0$.

\subsection*{p4mg}
\begin{figure}[btbp]
  \centering
  \begin{minipage}[t]{0.45\textwidth}   
  \centering
  \includegraphics[height=5cm]{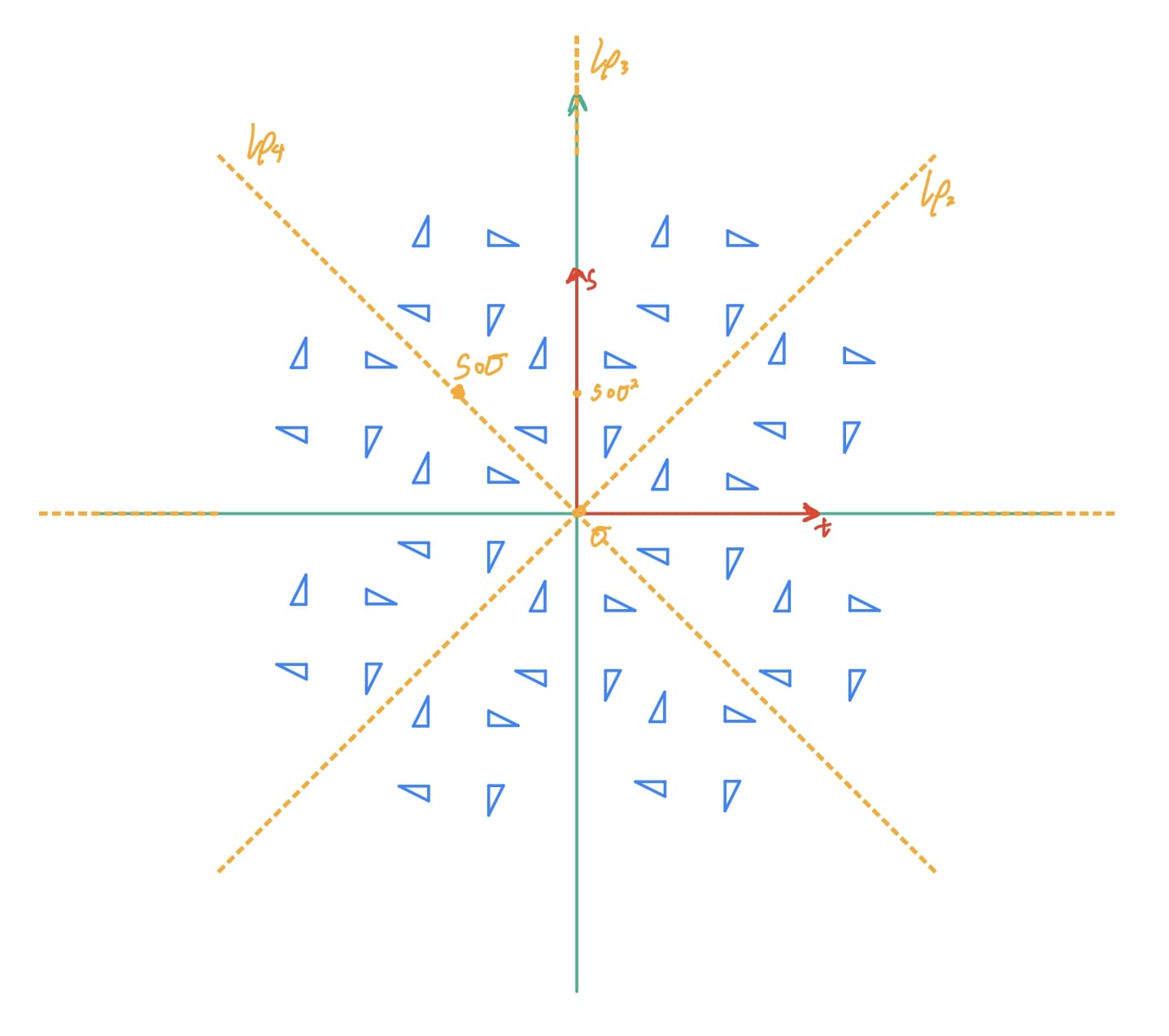}
  \caption{p4mg}  
  \end{minipage}
  \begin{minipage}[t]{0.45\textwidth}  
  \centering  
  \includegraphics[height=5cm]{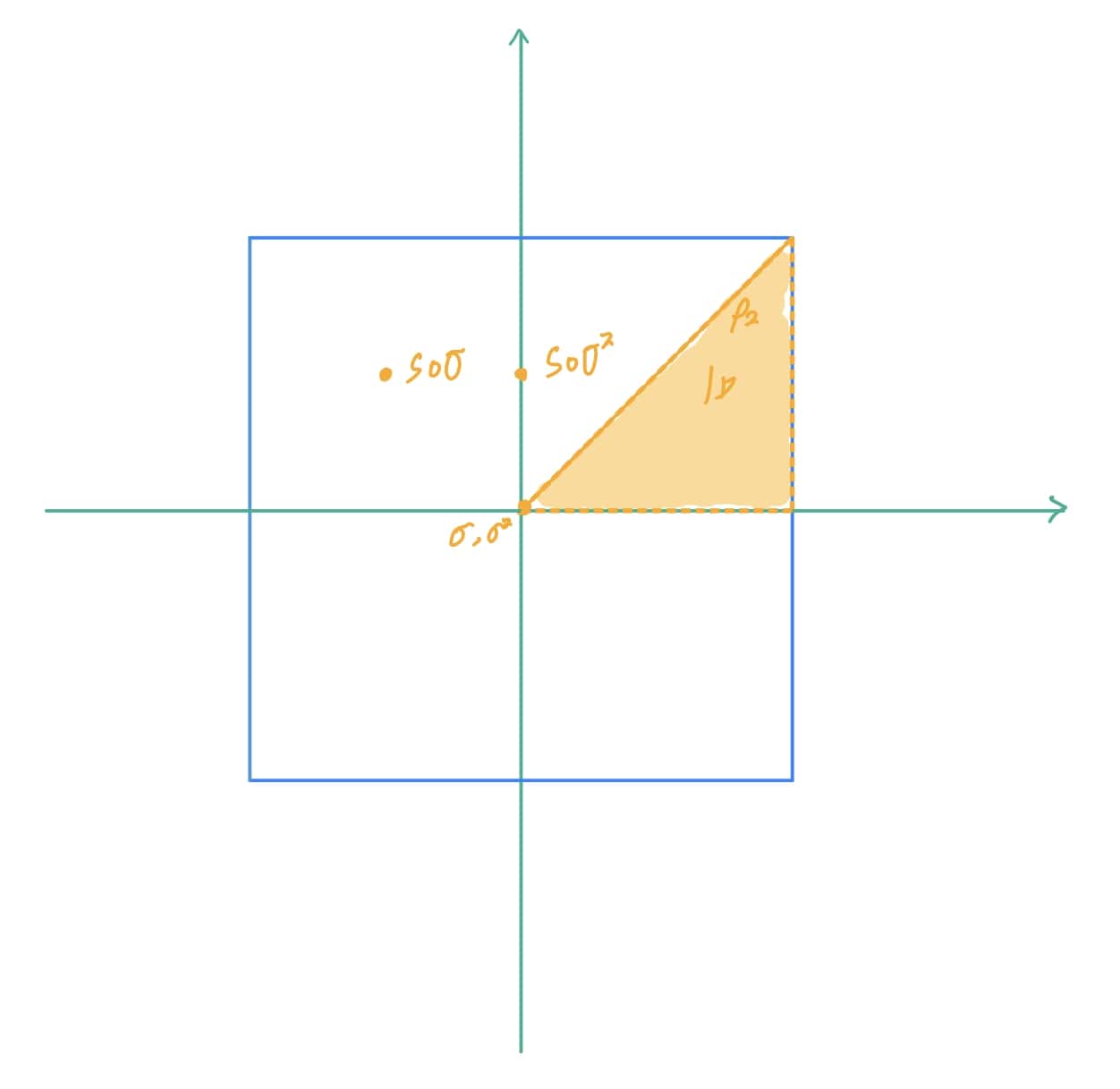}  
  \caption{Fundamental Domain of p4mg}  
  \end{minipage}  
  \end{figure}

We compute $K^{G}_{\bullet}(\underline{E}G)\cong K^{D}_{\bullet}(\mathbb{T}^{2})$.

1. Let $t=(t_{1},0)$ and  $s=(0,s_{2})$ form a basis for $A$ ($t_{1},s_{2}\neq0$ ). 

Point group $D=\{id_{D},\sigma,\sigma^{2},\sigma^{3},\rho_{1},\rho_{2},\rho_{3},\rho_{4}\}$ and the pullback of the genetors:
\begin{align*}
  \gamma(\sigma)=((0,0),M_{\sigma}),\quad \gamma(\rho_{1})=((\frac{1}{2},\frac{1}{2}),M_{\rho_{1}}), \\
M_{\sigma}=\begin{bmatrix}
0 & 1\\ 
-1 & 0
\end{bmatrix}, \quad
M_{\rho_{1}}=\begin{bmatrix}
  1 & 0\\ 
  0 & -1
  \end{bmatrix} .
\end{align*}

2. $FC(G)=FC(D)=\{[id_{D}],[\sigma],[\sigma^{2}],[\rho_{2}],[s\circ\sigma],[s\circ \sigma^{2}]\}$, where
$$
\begin{aligned}
s \circ \sigma =((0,1),M_{\sigma})=((-\frac{1}{2},\frac{1}{2})+(\frac{1}{2},-\frac{1}{2})M_{\sigma},M_{\sigma}),\\
s \circ \sigma^{2} =((0,1),M_{\sigma})=((0,\frac{1}{2})+(0,-\frac{1}{2})M_{\sigma^{2}},M_{\sigma^{2}}).
\end{aligned}
$$

3. See Table \ref{p4mg}.
\begin{table}[H]
  \footnotesize
  \centering%
  \begin{tabu} to 0.95\textwidth{X[c]X[c]X[c]X[c]}
  \hline
  representative         & identity      & $\rho_{2}$ & $\sigma$   \\ \hline
  $X^{g}$                & $\mathbb{T}^{2}$ & $(x,x)$ &$(0,0)$\\
  $C_{D}(g)$       & $D$    & $t+s$, $\rho_{2}$, $\rho_{4}$ & $\sigma$, $\sigma^{2}$  \\
  $X^{g}/C_{D}(g)$ & Fundamental Domain & $\{(x,x)|0\leqslant x\leqslant 1\}$&$(0,0)$\\
  \tiny{$\bigoplus \limits_{k \in even} H_{k}(X^{g}/C_{D}(g))$} & $\mathbb{Z}$ & $\mathbb{Z}$ & $\mathbb{Z}$  \\
  \tiny{$\bigoplus \limits_{k \in odd} H_{k}(X^{g}/C_{D}(g))$} &
  0&0&0\\ \hline \hline
  $\sigma^{2}$ & $s\circ\sigma$ &$s\circ\sigma^{2}$ & \\
  $(0,0)$ & $(-\frac{1}{2},\frac{1}{2})$ & $(0,\frac{1}{2})$ & \\
  $\sigma$, $\sigma^{2}$ & $s\circ\sigma$ & $s\circ\sigma^{2}$ & \\
  $(0,0)$ & $(-\frac{1}{2},\frac{1}{2})$ & $(0,\frac{1}{2})$ & \\
  $\mathbb{Z}$ & $\mathbb{Z}$ & $\mathbb{Z}$ & \\
  0&0&0& \\ \hline 
  \end{tabu}%
  \caption{The computation about p4mg}
  \label{p4mg}
  \end{table}

4. $K_{0}(p4mg)=\mathbb{Z}^{6}, K_{1}(p4mg)=0$.

\subsection*{p6mm}
\begin{figure}[btbp]
  \centering
  \begin{minipage}[t]{0.45\textwidth}   
  \centering
  \includegraphics[height=5cm]{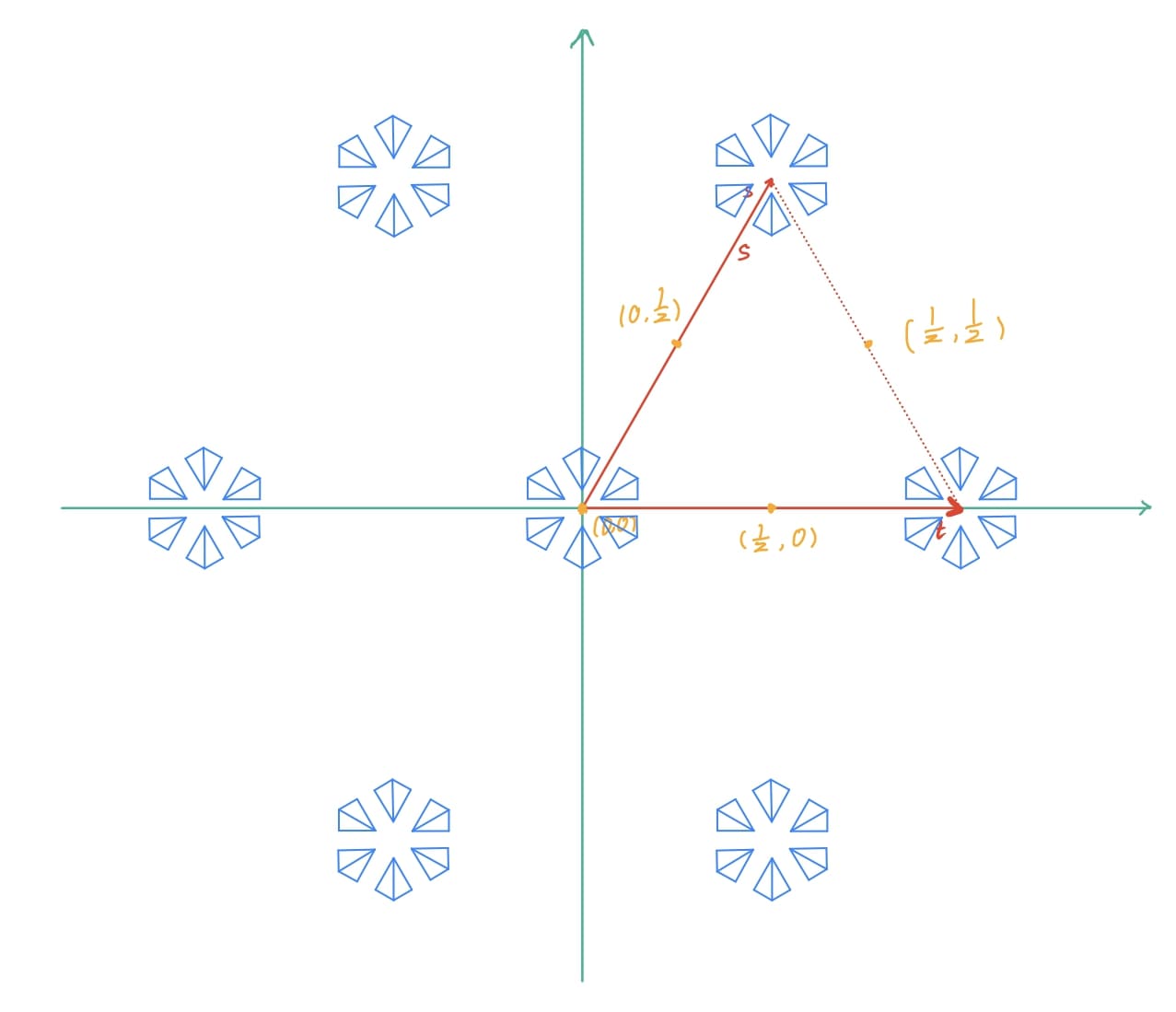}
  \caption{p6mm}  
  \end{minipage}
  \begin{minipage}[t]{0.45\textwidth}  
  \centering  
  \includegraphics[height=5cm]{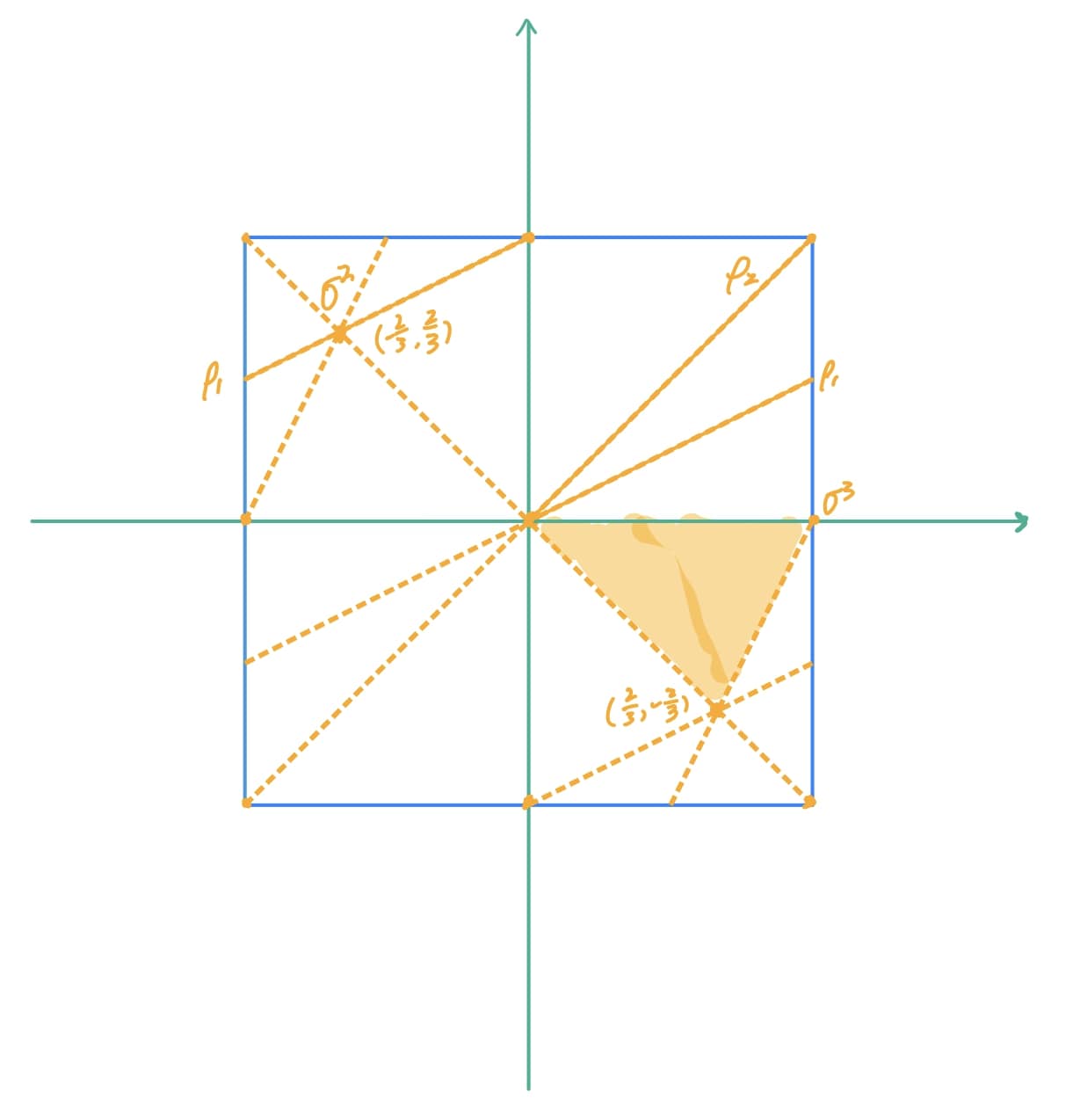}  
  \caption{Fundamental Domain of p6mm}  
  \end{minipage}  
  \end{figure}

We compute $K^{G}_{\bullet}(\underline{E}G)\cong K^{D}_{\bullet}(\mathbb{T}^{2})$.

1. Let $t=(t_{1},0)$ and  $s=(\frac{1}{2}t_{1},\frac{\sqrt{3}}{2}s_{2})$ form a basis for $A$ ($t_{1},s_{2}\neq0$ ). 

Point group $D=\{id_{D},\sigma,\sigma^{2},\sigma^{3},\sigma^{4},\sigma^{5},\rho_{1},\rho_{2},\rho_{3},\rho_{4},\rho_{5},\rho_{6}\}$ and the pullback of the genetors:
\begin{align*}
  \gamma(\sigma)=((0,0),M_{\sigma}),\quad \gamma(\rho_{1})=((0,0),M_{\rho_{1}}), \\
M_{\sigma}=\begin{bmatrix}
0 & 1\\ 
-1 & 1
\end{bmatrix}, \quad
M_{\rho_{1}}=\begin{bmatrix}
  1 & 0\\ 
  1 & -1
  \end{bmatrix} .
\end{align*}

2. $FC(G)=FC(D)=\{[id_{D}],[\sigma],[\sigma^{2}],[\sigma^{3}],[\rho_{1}],[\rho_{2}]\}$.

3. See Table \ref{p6mm}.
\begin{table}[H]
  \footnotesize
  \centering%
  \begin{tabu} to 0.95\textwidth{X[c]X[c]X[c]X[c]}
  \hline
  representative         & identity      & $\sigma$ & $\sigma^{2}$   \\ \hline
  $X^{g}$                & $\mathbb{T}^{2}$ & $(0,0)$ &$(0,0)$, $(\frac{2}{3},-\frac{2}{3})$, $(-\frac{2}{3},\frac{2}{3})$\\
  $C_{D}(g)$       & $D$    & $id_{D}$, $\sigma$, $\sigma^{2}$, $\sigma^{3}$, $\sigma^{4}$, $\sigma^{5}$ & $id_{D}$, $\sigma$, $\sigma^{2}$, $\sigma^{3}$, $\sigma^{4}$, $\sigma^{5}$  \\
  $X^{g}/C_{D}(g)$ & Fundamental Domain & $(0,0)$&$(0,0)$, $(\frac{2}{3},\frac{2}{3})$\\
  \tiny{$\bigoplus \limits_{k \in even} H_{k}(X^{g}/C_{D}(g))$} & $\mathbb{Z}$ & $\mathbb{Z}$ & $\mathbb{Z}^{2}$  \\
  \tiny{$\bigoplus \limits_{k \in odd} H_{k}(X^{g}/C_{D}(g))$} &
  0&0&0\\ \hline \hline
  $\sigma^{3}$ & $\rho_{1}$ &$\rho_{2}$ & \\
  $(\pm 1,\pm 1)$, $(0,0)$, $(0,\pm 1)$, $(\pm 1,0)$ & $(x,\frac{x}{2})$, $(x,\frac{x}{2}\pm 1)$  & $(x,x)$ & \\
  $id_{D}$, $\sigma$, $\sigma^{2}$, $\sigma^{3}$, $\sigma^{4}$, $\sigma^{5}$ & $id_{D}$, $\rho_{1}$, $\sigma^{3}$ & $id_{D}$, $\rho_{2}$, $\sigma^{3}$ & \\
  $(0,0)$, $(1,0)$&$\{(x,\frac{x}{2})|0\leqslant x\leqslant 1\}$, $\{(x,\frac{x}{2}+1)|-1\leqslant x\leqslant 0\}$ &$\{(x,x)|0\leqslant x\leqslant 1\}$& \\
  $\mathbb{Z}^{2}$ & $\mathbb{Z}$ & $\mathbb{Z}$ & \\
  0&0&0& \\ \hline 
  \end{tabu}%
  \caption{The computation about p6mm}
  \label{p6mm}
  \end{table}

4. $K_{0}(p6mm)=\mathbb{Z}^{8}, K_{1}(p6mm)=0$.

\bibliographystyle{unsrt}
\bibliography{ref}

\end{document}